\documentclass[preprint,8pt]{elsarticle}

\usepackage{amssymb}

\usepackage{typearea}
\typearea{1}

\usepackage{syntonly}

\usepackage[paperwidth=16.5cm,
            paperheight=24cm,
	    		text={12.5cm,18.5cm},headheight=1cm]{geometry}    
\usepackage{latexsym}                      
\usepackage[english]{babel}                
\usepackage{subfigure}                     
\usepackage{graphics}
\usepackage{graphicx}

\usepackage{stmaryrd}
\usepackage{booktabs}                      
\usepackage{fancyhdr}                      
\usepackage{prettyref}                     
\usepackage{overpic}                       
\usepackage{amssymb}                       
\usepackage{amsmath}                       
\usepackage{amsthm}                        
\usepackage{amstext}                       
\usepackage{array}                         
\usepackage{natbib} 
\usepackage{ifthen}                        
\usepackage{ifpdf}                         
\usepackage[usenames]{color}               
\usepackage{rotating}                      
\usepackage{hyphenat}                      
\usepackage{makeidx}                       

\usepackage{fnbreak}                       
\usepackage{xspace}                        
\usepackage{multirow}
\usepackage{fancyvrb} 
\usepackage{nomencl}                     
\usepackage{url}                           
\usepackage[active]{srcltx}                

\usepackage{tikz}
\usetikzlibrary{shapes,arrows}

\usepackage[titletoc,toc,title]{appendix}

\ifpdf
	\usepackage[pdftex,
		    pagebackref=true,
		    pdfborder=000,
		    bookmarksnumbered=true,
        plainpages=false,
		    pdfpagelabels=true,
        pdftitle={TBC}, 
		    pdfauthor={Jes\'us Garicano-Mena},
		    pdfsubject={Research article},
		    pdfkeywords={computational fluid dynamics, heat transfer, residual distribution schemes, unstructured grids, hypersonic flow,thermo-chemical non-equilibrium}
		    ]{hyperref}
\hypersetup{colorlinks=false}
\else
	\usepackage[hypertex,pagebackref=true]{hyperref}
\fi

\usepackage[normalem]{ulem} 
\usepackage{adjustbox} 

\newcommand{\be}{\begin{equation}}
\newcommand{\ee}{\end{equation}}
\newcommand{\bi}{\begin{itemize}}
\newcommand{\ei}{\end{itemize}}

\usepackage{framed}                      %

\usepackage[english]{babel}
\usepackage{subfigure}
\newlength{\figwidth}
\newlength{\minifigwidth}
\usepackage{float}

\floatstyle{boxed} 
\newfloat{algorithm}{htbp}{lop}
\floatname{algorithm}{Algorithm}

\usepackage{xfrac}

\newcommand{\beq}{\begin{equation}}
\newcommand{\eeq}{\end{equation}}

\usepackage{amsmath}
\usepackage{mathtools}
\usepackage{stmaryrd}
\usepackage{lineno}
\usepackage{xcolor}
\modulolinenumbers[5]
\journal{Journal of Computational Physics}

\begin{document}

\begin{frontmatter}



\title{On the role of numerical dissipation in stabilising under--resolved turbulent simulations using discontinuous Galerkin methods}



%
%
%



		\author[upmaddress]{Juan Manzanero\corref{mycorrespondingauthor}}
		\cortext[mycorrespondingauthor]{Corresponding author}
		\ead{juan.manzanero@upm.es}		
		\author[upmaddress,ccsaddress]{Esteban Ferrer}
		\author[upmaddress,ccsaddress]{Gonzalo Rubio}
		\author[upmaddress,ccsaddress]{Eusebio Valero}				

		\address[upmaddress]{ETSIAE-UPM - School of Aeronautics, Universidad Polit\'ecnica de Madrid,\\ Plaza Cardenal Cisneros 3, E-28040 Madrid, Spain. }
\address[ccsaddress]{CCS-UPM - Center for Computational Simulation, Universidad Polit\'ecnica de Madrid,\\ Boadilla del Monte, 28660 Madrid, Spain.}

\begin{abstract}
We analyse numerical errors (dissipation and dispersion) introduced by the discretisation of inviscid and viscous terms in energy stable discontinuous Galerkin methods. First, we analyse these methods using a linear von Neumann analysis (for a linear advection--diffusion equation) to characterise their properties in wave--number space. Second, we validate these observations using the 3D Taylor--Green Vortex Navier--Stokes problem to assess transitional/turbulent flows.

We show that the dissipation introduced by upwind Riemann solvers affects primarily high wave--numbers. 
 This dissipation may be increased, through a penalty parameter, until a critical value. However, further augmentation of this parameter leads to a decrease of dissipation, reaching zero for very large values.
Regarding the dissipation introduced by second order derivatives, we show that this dissipation acts at low and medium wave--numbers (lower wave--numbers compared to upwind Riemann solvers).

 In addition, we analyse the Spectral Vanishing Viscosity (SVV) technique, previously used in continuous discretisations (e.g. Fourier), to find that with an appropriate kernel (which damps selected modes) it is possible to control the amount of dissipation introduced in the low and medium wave--number range.

%

%

%
%
Combining these ideas, we finally propose a DG-SVV approach that uses a Smagorinsky model to compute the numerical viscosity. This DG-SVV approach is tested in an isotropic laminar/turbulent under--resolved scenario. Combining the SVV technique with a low dissipation Riemann solver, we obtain a scheme capable of maintaining low dissipation levels for laminar flows, whilst providing the correct dissipation for all wave--number ranges in turbulent regimes.
\end{abstract}

\end{frontmatter}

\tableofcontents
\section{Introduction}\label{sec:introduction}
High order discontinuous Galerkin (DG) methods have been adopted by academia and research centers as an alternative to classic numerical schemes (e.g. Finite Differences, Finite Volumes or Finite Elements). During recent years, DG methods have been adapted to solve increasingly complex physics; including incompressible  and compressible flow problems \cite{BR1,2017:Ferrer,Kompenhans2016216,Kompenhans201636,Ferrer2010,Ferrer2012,
Ferrer_CiCP,Fraysse2016805}.
This popularity may be attributed to two characteristics. First, DG methods provide high accuracy even for unstructured distorted meshes \cite{2008:Hesthaven}, a property difficult to retain when using classic methods (e.g. Finite Differences or Finite Volumes). This property is a result of the compactness and local character of the scheme, enabling high order accuracy using compact stencils.
Second, DG methods have shown to be more robust that their high order continuous Galerkin relatives \cite{2003:Kirby}.
 This last beneficial property  has been often attributed to the use of upwind Riemann solvers, which add controlled dissipation, enhancing robustness.
In recent years, the DG community has exploited the increased robustness provided by Riemann solvers to simulate  under--resolved turbulent flows, e.g. \cite{2013:GassnerUnder,2017:Moura_2}. The idea behind pushing robustness towards higher Reynolds numbers is to rely on the localised dissipation provided by upwind Riemann solvers to dissipate small flow structures, which cannot be resolved on coarse meshes. Methods that rely on the numerics to provide mechanisms for turbulent dissipation are typically known as implicit Large Eddy Simulation (iLES) methods. The term implicit evidences that
numerical errors (in particular dissipation) are in charge of dissipating under--resolved flow structures. An alternative to iLES is provided by explicit LES methods, where flow dissipation at small scales relies on physical arguments and modified flow equations, e.g. \cite{2000:Lenormand,2001:Hughes}.
 When solving under--resolved turbulent flows using iLES methods, it is necessary to understand and control numerical errors, and specially numerical dissipation introduced by the scheme, which replaces explicit subgrid--scale models.
%

 In DG methods, there are several alternatives to provide numerical dissipation.
 The most popular choice is to include upwind Riemann solvers, which arise naturally from the integration by parts of the non--linear fluxes and the existence of inter--element discontinuities in DG. These fluxes introduce local dissipation, which scales with the size of the discontinuities in the numerical solution. In under--resolved flows, the size of such discontinuities increases and it has been argued (see for example \cite{2017:Ferrer}) that fluxes based on discontinuities may act as an appropriate stabilising mechanism for under--resolved turbulent simulations.
 In addition to upwind Riemann solvers, we consider artificial viscosity methods \cite{1950:Neum} and the Spectral Vanishing Viscosity (SVV) technique \cite{2006:Kirby} to include dissipation, when the dissipation resulting from the Riemann solver is insufficient. Artificial viscosity methods add a controlled amount of viscosity to the governing equations in the vicinity of strong gradients. As far as the SVV technique is concerned, it was first introduced for  Fourier and continuous Galerkin discretisations \cite{1989:Tadmor} to regularise the solution (i.e. avoid oscillatory phenomena)
  in the inviscid Burgers equation, and later in the Navier--Stokes equations \cite{2000:Karamanos}. SVV is similar in spirit to including artificial viscosity, and provides additional dissipation (only at high wave-numbers) that enhances stability, vanishes in the laminar limit and provides spectral convergence in high--order discretisations.
Let us note that other techniques to introduce localised dissipation in continuous discretisations exist but are not considered in this work, e.g. SUPG  stabilisation.

 A substantial amount of work \cite{2013:Hillewaert,2017:Ferrer,2017:Flad,2017:Fehn,2017:DeLaLlave} has focused on understanding the stabilising effect of upwind Riemann solvers on under--resolved turbulent flows.  However, limited efforts have been devoted to understand the combined effect of the Riemann solvers and discretised viscous terms.
 In this text, we analyse the individual contributions of each term but also the combined effect.
For all dissipative mechanisms, we first analyse their numerical properties in wave--number space, using von Neumann analysis on a linear advection--diffusion equation. Similar work has been presented for linear advection with constant coefficients \cite{2011:Gassner,2015:Moura}  and for non--constant coefficients
\cite{2016:Manzanero}. Second, we correlate these findings with results for DNS and iLES simulations of the 3D Taylor--Green Vortex (TGV) Navier--Stokes problem with transitional/turbulent flow \cite{2013:GassnerUnder,2017:Moura_2}. Even though von Neumann analyses are restricted to the constant advection--diffusion equation, we will confirm, with the help of numerical experiments, that von Neumann results are consistent with observations of under--resolved Navier--Stokes turbulent flows.
We show that the dissipation introduced by upwind Riemann solvers affects high wave--numbers, whilst discrete second order derivatives provide dissipation at low and medium wave--numbers. Additionally, the dissipation introduced by SVV operators helps to control the amount of dissipation introduced in the low and medium wave--number range, and provides a suitable mechanism to develop new models, as the one proposed in the final section of this manuscript.
To perform these studies, we start from a baseline scheme without dissipation (i.e. an energy conserving scheme), and include numerical dissipation through  the different stabilisation techniques. In the linear advection equation with constant coefficient it suffices to consider central fluxes (see \cite{Kopriva2,2017:Manzanero}), but in the non--linear Navier--Stokes equations, it is required to use split formulations \cite{2016:gassner}, to minimise aliasing errors.

Having quantified numerical errors for the above dissipative mechanisms, in a final section, we combine upwind Riemann solvers and spectral vanishing viscosity. Following the suggestion of Karamanos and Karniadakis for continuous Galerkin methods \cite{2000:Karamanos,2002:Kirby}, we modify the classic SVV technique using a Smagorinsky model to adjust the amount of dissipation introduced.
This new proposed DG model is capable of maintaining low dissipation levels in laminar flows, whilst modelling small eddies and providing correct dissipation for all wave--number ranges in turbulent regimes.


The rest of this paper is organised as follows: the description of the dissipative mechanisms and their inclusion in the Navier--Stokes equations and the 1D advection--diffusion equation, are included in Section \ref{sec:problem_description}.
In Section \ref{sec:VN_and_TGV} we study three dissipation techniques: in Section \ref{subsec:results_upwindRS} we analyse upwind Riemann solvers, in Section \ref{subsec:results_LES} we investigate the discretisation of viscous terms , and in Section \ref{subsec:results_SVV} we assess artificial viscosity (LES models and spectral vanishing viscosity) stabilisation. Finally, we propose and test a Smagorinsky--SVV DG discretisation to simulate under--resolved turbulent/transitional flows in Section \ref{sec:LESSVV_hybrid}.

\section{Methodology} \label{sec:problem_description}


We first present, in Section \ref{subsec:advectionAndVonNeumann},  the 1D advection--diffusion equation together with von Neumann stability analysis. In Section \ref{subsec:3DNSE} the discretisation of the 3D compressible Navier--Stokes equations is sketched. Then in Section \ref{subsec:AboutECSchemes} we briefly introduce energy conserving schemes. Finally in Section \ref{subsec:DissipationTechniques} we detail how the dissipation mechanisms are included in both the 3D compressible Navier--Stokes equations and the 1D advection--diffusion equation: upwind Riemann solvers in Section \ref{subsubsec:upwindRiemannTheory}, discretisation of viscous terms in Section \ref{subsubsec:theory_LES}, and artificial viscosity (LES models and spectral vanishing viscosity) in Section \ref{subsubsec:theory_SVV}.

\subsection{1D advection--diffusion equation and von Neumann Analysis}\label{subsec:advectionAndVonNeumann}

The 1D advection--diffusion equation reads
\begin{equation}\label{eq:AdvEq}
u_t+au_x=(\mu u_x)_x,
\end{equation}
where $a$ is a constant advection speed and $\mu u_x$ (being $\mu$ the viscosity) is the viscous flux. Even though the viscosity, $\mu$, will be kept constant in what follows, we preferred to consider the divergence form for the viscous flux, $\mu u_x$, for parallelisms with the viscous fluxes in the Navier--Stokes equations (see next section).
Following \cite{2009:Kopriva}, a weak form can be constructed by multiplying \eqref{eq:AdvEq} by a smooth test function $\phi$ and integrating by parts in each element, $el$:

\begin{equation}
\int_{{el}}\phi u_t  + \underbrace{\phi af^\star}_\text{Flux dissipation} \biggr|_{\partial el} -\int_{{el}}\phi_x au  = \int_{{el}} \phi (\underbrace{(\mu u_x)_x}_\text{Viscous term} + \underbrace{(f_{diss})_x}_\text{Artificial dissipation}),
\label{eq:DGScheme}
\end{equation}
where $\partial el$ represents the element boundaries. Notice that equation \eqref{eq:DGScheme} has been augmented with an artificial dissipative term $(f_{diss})_x$.
The numerical flux function $f^\star$ arises as a result of the discontinuities of the numerical solution across the interfaces.
This flux can provide dissipation, depending on the discretisation used. For example, upwinding provides dissipation, whilst central fluxes do not. Further details are included in following sections.

The 1D advection--diffusion equation is used for von Neumann analyses to quantify dispersion and dissipation errors introduced by the numerical discretisation. Our analysis provides insight into the effect of the different dissipation techniques described in Section \ref{subsec:DissipationTechniques}.
In particular, we consider the discretisation of the linear advection equation \eqref{eq:AdvEq} by means of the discontinuous Galerkin scheme, equation \eqref{eq:DGScheme}, in an uniform grid with element spacing $h$.
For the sake of simplicity, only the advective terms of equation \eqref{eq:AdvEq} (LHS terms) are considered in the description of von Neumann analysis, however the diffusive terms (RHS terms in equation \eqref{eq:AdvEq}) are also taken into account to obtain the results. Details on the derivation of von Neumann analysis for an advection--diffusion equation may be found in \cite{2016:Moura,2016:Manzanero}.
To perform von Neumann analysis we introduce an exponential wave solution, with spatial wave--number $k$,

\begin{equation}
u(x,t) = e^{ikx-i\omega t} = u_0(x)e^{-i\omega t},\label{eq:InitialConditionVN}
\end{equation}
such that the temporal frequency, $\omega$, yields an eigenfunction of the original PDE \eqref{eq:AdvEq} (i.e. it satisfies $\omega=ak$). Because of the linear relation between spatial and temporal frequencies, we will only consider a constant coefficient $a=1$. The choice of an exponential initial condition \eqref{eq:InitialConditionVN} is essential, because it allows to relate the solution between elements by a spatial phase shift:

\begin{equation}
\{\underline{\boldsymbol{u}}^{el-n}\} = e^{-iknh}\{\underline{\boldsymbol{u}}^{el}\}. \label{eq:relationship_disp_diff}
\end{equation}
where $\{\underline{\boldsymbol{u}}^{el}\}$ is a vector containing the nodal degrees of freedom of the $el-$th element. Introducing \eqref{eq:relationship_disp_diff} in the discrete version of \eqref{eq:DGScheme}, we obtain a linear ordinary differential equation system for each individual element:

\begin{equation}
\frac{h}{2}\frac{d}{dt}\{\underline{\boldsymbol{u}}\} = \bigl(e^{-ikh}[\boldsymbol{L}]+[\boldsymbol{C}]+e^{ikh}[\boldsymbol{R}]\bigr)\{\underline{\boldsymbol{
u}}\} = [\boldsymbol{M}(kh)]\{\underline{\boldsymbol{u}}\}.
\label{eq:vonNeumann}
\end{equation}
The precise expression of matrices $[\boldsymbol{L}]$, $[\boldsymbol{C}]$, and $[\boldsymbol{R}]$ can be found in \cite{2016:Manzanero}. The general solution of \eqref{eq:vonNeumann} is linearly spanned by the $N+1$ modes of the eigenvalue problem (note that, for simplicity, the index $el$ has been dropped):

\begin{equation}
-i\omega \frac{h}{2}\{\underline{\boldsymbol{v}}\} = [\boldsymbol{M}(kh)]\{\underline{\boldsymbol{v}}\}, \quad  u = \sum_{m=0}^N A_m \{\underline{\boldsymbol{v}_m}\}e^{-i\omega_m t},
\label{eq:eigenvalueproblem}
\end{equation}
where $\{\boldsymbol{v}_m\}$ are the unitary eigenvectors, and the amplitudes $A_m$ are constants computed to recover the initial condition in $t=0$:

\begin{equation}
\{\underline{\boldsymbol{u}}(t=0)\} = \sum_{m=0}^N A_m \{\underline{\boldsymbol{v}}_m\}e^{-i\omega_m t}\biggr|_{t=0} =  \sum_{m=0}^N A_m \{\underline{\boldsymbol{v}}_m\} = \{\underline{\boldsymbol{u}}_0\}.
\label{eq:initial_condition_modes}
\end{equation}
The solution structure \eqref{eq:initial_condition_modes} allows us to classify three different numerical error sources, which were already detailed in \cite{2015:Moura,2016:Manzanero}. Specifically, only the so-called \textit{primary mode} (p) propagates with the physical wave-speed and damping (e.g. $\omega_p=0$ when $k=0$). Hence, we rewrite the solution isolating the primary mode contribution from the remaining modes (called \textit{secondary modes}, $m\neq p$):

\begin{equation}
\{\underline{\boldsymbol{u}}\} = A_p \{\underline{\boldsymbol{v}}_p\}e^{-i\omega_p t} + \sum_{\substack{m=0 \\ m\neq p}}^N A_m \{\underline{\boldsymbol{v}}_m\}e^{-i\omega_m t}.\label{eq:modes_primary}
\end{equation}
Next, the initial condition \eqref{eq:initial_condition_modes} is also separated in the contribution of primary and secondary modes, and introduced in \eqref{eq:modes_primary}. As a result, we can consider the numerical solution as the primary mode propagating the initial condition, $ \{\underline{\boldsymbol{u}}_0\}e^{-i\omega_p t}$, and the non--physical errors as secondary modes, $ \{\Delta \underline{\boldsymbol{u}}(t)\}$:

\begin{equation}
\begin{split}
\{\underline{\boldsymbol{u}}\} &=  \{\underline{\boldsymbol{u}}_0\}e^{-i\omega_p t} + \sum_{\substack{m=0 \\ m\neq p}}^N A_m \{\underline{\boldsymbol{v}}_m\}\bigl(e^{-i\omega_m t} - e^{-i\omega_p t}\bigr)\\
 &=\{\underline{\boldsymbol{u}}_0\}e^{-i\omega_p t} + \{\Delta \underline{\boldsymbol{u}}(t)\}.
\end{split}
\label{eq:secondary_modes_error}
\end{equation}

Besides, the numerical propagation speed experienced by the initial condition, $\{\underline{\boldsymbol{u}}_0\}$, (i.e. the primary mode frequency, $\omega_p$) will differ from that dictated by the analytical PDE ($\omega=ak$). Precisely, the difference between its real part and the theoretical travelling speed, $ak$, will yield a dispersion error (i.e. error in the propagation speed), whilst its imaginary part, which is generally non--zero, entails numerical dissipation.
Following \cite{2015:Moura}, we define the non--dimensional wave--number, $\hat{k}$, as:

\begin{equation}
\hat{k} = \frac{kh}{N+1},
\end{equation}
such that we will distinguish between low wave--numbers $(\hat{k}<\pi/2)$, medium wave--numbers $(\hat{k}\sim\pi/2)$ and high wave--numbers $(\hat{k}                                                                                                                                                                                                                                                                                                                                                         >\pi/2)$.

 Further details on the discretisation and analysis of the 1D advection equation and von Neumann analysis for advection equations with non--constant coefficients may be found in other works by the authors in \cite{2016:Manzanero,2017:Manzanero,2018:Manzanero}.

\subsection{3D Navier--Stokes equations}\label{subsec:3DNSE}

The 3D Navier--Stokes equations can be compactly written as:
%
\begin{equation}\label{eq::NS1}
\vec{u}_t   +\nabla\cdot\vec{\boldsymbol{F}}_e= \nabla\cdot\vec{\boldsymbol{F}}_v,
\end{equation}
where $\vec{u}$ is the vector of conservative variables $\vec{u} = [ \rho , \rho v_1 , \rho v_2 , \rho v_3 , \rho e]^T$. Details on the specific formulations retained for inviscid and viscous fluxes $\vec{\boldsymbol{F}}_e$ and $\vec{\boldsymbol{F}}_v$ can be found in \ref{appendix::A} of this text.

To derive discontinuous Galerkin schemes, we consider \eqref{eq::NS1} for one mesh element $el$, multiply by a locally smooth test function $\phi_j$, for $0\leq j\leq N$, where $N$ is the high order polynomial degree, and integrate on $el$:
\begin{equation}\label{eq::NS2}
\int_{el}\vec{u}_t\phi_j+\int_{el} \nabla\cdot\vec{\boldsymbol{F}}_e\phi_j  =\int_{el} \nabla\cdot\vec{\boldsymbol{F}}_v\phi_j.
\end{equation}
We can now integrate by parts the inviscid fluxes, $\boldsymbol{F}_e$, integral to obtain a local weak form of the equations (one per mesh element):
\begin{equation}\label{eq::NS3}
\int_{el}\vec{u}_t\phi_j +  \int_{\partial el} \vec{\boldsymbol{F}}_e\cdot\mathbf{n}\phi_j-\int_{el} \vec{\boldsymbol{F}}_e\cdot\nabla\phi_j
=\int_{el} \nabla\cdot\vec{\boldsymbol{F}}_v\phi_j,
\end{equation}
where $\mathbf{n}$ is the normal vector at element boundaries ${\partial el}$. 
We replace discontinuous fluxes at inter--element faces by a numerical inviscid flux,  ${\vec{\boldsymbol{F}}^*}_e$, to obtain a weak form for the equations for each element, 

\begin{equation}\label{eq::NS4}
\int_{el}\vec{u}_t\cdot\phi_j + \int_{\partial el} {\vec{\boldsymbol{F}}^*}_e\cdot\mathbf{n}\phi_j-\int_{el} \vec{\boldsymbol{F}}_e\cdot\nabla\phi_j
=\int_{el} \nabla\cdot\vec{\boldsymbol{F}}_v\phi_j,
\end{equation}
where, we have omitted the fluxes at external boundaries, for simplicity. This set of equations for each element is coupled through the inviscid fluxes ${\vec{\boldsymbol{F}}^*}_e$ and governs flow behaviour. Finally, equation \eqref{eq::NS4} may be augmented with an artificial dissipation term $\nabla\cdot\vec{\boldsymbol{F}}_{diss}$: 

\begin{equation}\label{eq::NS5}
\int_{el}\vec{u}_t\cdot\phi_j + \int_{\partial el} \underbrace{{\vec{\boldsymbol{F}}^*}_e\cdot\mathbf{n}}_\text{Riemann solver}\phi_j-\int_{el} \vec{\boldsymbol{F}}_e\cdot\nabla\phi_j
=\int_{el} (\underbrace{\nabla\cdot\vec{\boldsymbol{F}}_v}_\text{Viscous term}+\underbrace{\nabla\cdot\vec{\boldsymbol{F}}_{diss}}_\text{Artificial dissipation})\cdot\phi_j,
\end{equation}
where the various terms that can be discretised to control dissipation in the numerical scheme have been underlined. Details on the various forms of such terms are included in following sections. Additionally, Table \ref{my-label} summarises the various dissipative mechanisms and details their effect wave--numbers space.

\begin{table}[!]
	\centering
	\caption{Summary of dissipative mechanisms: $\hat{k}$ denotes the non--dimensional wave--number, $\hat{k}=kh/(N+1)$. }
	\label{my-label}
\small
	\begin{tabular}{l|c|c|c}
		& Discretisation                                         & Section & Numerical dissipation range                \\
\hline
		upwind   & Roe, low   &   \ref{subsubsec:upwindRiemannTheory}& High wave--numbers        \\
		 Riemann solver &  dissipation Roe \cite{2016:Obwald}   & \ref{subsec:results_upwindRS}  & $(\hat{k}>\pi/2)$         \\
\hline
		Viscous  & Bassi-Rebay (BR1)   & \ref{subsubsec:theory_LES}&Low $(\hat{k}<\pi/2)$ and  medium\\
		term &   &  \ref{subsec:results_LES} &     $(\hat{k}\sim\pi/2)$     wave--numbers \\
\hline
		Artificial & LES models; SVV  & \ref{subsubsec:theory_SVV}& Low $(\hat{k}<\pi/2)$ and  medium      \\
		 dissipation   &   & \ref{subsec:results_SVV}  &$(\hat{k}\sim\pi/2)$     wave--numbers \\
		    &   &   &(tuneable for SVV) \\		
	\end{tabular}
\normalsize
\end{table}

Further details regarding the nodal discontinuous Galerkin formulation used in this paper may be found in \cite{1999:Black} and in \cite{2009:Kopriva,2016:gassner} for extensions to 3D unstructured grids with curved elements.

\subsection{Briefly on energy conserving schemes} \label{subsec:AboutECSchemes}

Throughout the paper we use only energy conserving schemes, which are constructed to conserve discrete energy (assuming zero physical viscosity in \eqref{eq:AdvEq} or \eqref{eq::NS1}). Because the schemes are designed to remain stable and energy conserving, they do not require  numerical dissipation. Hence, including numerical dissipation to these schemes enables the analysis of dissipation techniques, which will not be masked by other dissipative mechanism. Overall, these schemes are useful to understand numerical dissipation requirements, and thus, to design robust dissipation techniques capable to achieve accurate solutions. Next, we introduce the dissipation techniques studied in this paper: upwind Riemann solvers, discretisation of viscous terms, and artificial viscosity (LES models and spectral vanishing viscosity).

In the linear advection equation, energy conserving schemes are achieved with an appropriate choice of the numerical Riemann solver. The numerical flux \eqref{eq:numericalflux} is designed such that the first part (central fluxes) balances volume terms to obtain an energy conserving scheme (for constant advection speeds). Hence, we will refer to the particular case with $\lambda=0$ as central fluxes.

The extension to the non--constant advection speed and non--linear problems (e.g. the Navier--Stokes equations) needs special treatment of volume integrals (to use split--forms) and requires the use of Gauss--Lobatto points to cancel out boundary terms using the summation--by--parts simultaneous--approximation--term property (SBP--SAT). The interested reader is referred to \cite{2013:Gassner,Kopriva2,2016:Manzanero,2016:gassner} for particular split--forms for non--constant advection, Burgers and Euler equations.

\subsection{Three dissipation mechanisms}\label{subsec:DissipationTechniques}
In this section, we introduce three dissipation mechanisms: upwind Riemann solvers in Section \ref{subsubsec:upwindRiemannTheory}, discretisation of viscous terms in Section \ref{subsubsec:theory_LES}, and artificial viscosity (LES models and spectral vanishing viscosity) in Section \ref{subsubsec:theory_SVV}. We correlate these mechanisms when included in the 1D advection--diffusion equation, to their counterparts in the  Navier-Stokes equations, such that in following sections, the information extracted from von Neumann analyses can be used to understand the behaviour when solving laminar/turbulent flows.

\subsubsection{First dissipation technique: upwind Riemann solvers}\label{subsubsec:upwindRiemannTheory}
 Energy conserving schemes present an intrinsic problem: the existence of non--dissipated, poorly approximated numerical modes.
 In this section we will add numerical dissipation through the inviscid numerical fluxes, $F^\star_e$, (see  \eqref{eq:numericalflux}). For the 1D advection equation with constant coefficient $a$, \eqref{eq:AdvEq}, the numerical flux is often defined as

\begin{equation}
f^\star(u_L,u_R) = a\frac{u_L+u_R}{2} - \frac{1}{2}\lambda |a|(u_R-u_L) = a\{u\} + \frac{1}{2}\lambda a \llbracket u\rrbracket\cdot\boldsymbol{n}_L,
\label{eq:numericalflux}
\end{equation}
where the following jump and average compact expressions have been used:

\begin{equation}
\{u\} = \frac{u_L + u_R}{2},~~ \llbracket u \rrbracket = u_L \boldsymbol{n}_L + u_R \boldsymbol{n}_R,
\end{equation}
and $\boldsymbol{n}_L$ and $\boldsymbol{n}_R$ are the outward pointing normal vectors of left and right elements respectively. We define the operator $\lambda\text{diss}(u,\phi)$ as the result of combining the weak integrals associated to the second part of the Riemann solver \eqref{eq:numericalflux} to a boundary shared by two elements:

\begin{equation}
\lambda\text{diss}(u,\phi) = -\frac{1}{2}|a|\lambda\llbracket u\rrbracket \cdot \boldsymbol{n}_L\phi_L -\frac{1}{2}|a|\lambda\llbracket u\rrbracket \cdot \boldsymbol{n}_R\phi_R = -\lambda
\frac{1}{2}|a|\llbracket u\rrbracket \llbracket \phi\rrbracket.
\label{eq:lambdastab}
\end{equation}
This contribution always provides numerical dissipation to the energy balance, obtained with $\phi=u$ and assuming $\lambda\geq 0$:

\begin{equation}
\lambda\text{diss}(u,\phi=u) = -\lambda\frac{1}{2}|a|\llbracket u\rrbracket^2 \leq 0.
\label{eq:lambdastabdiss}
\end{equation}
Additionally, note that $\lambda\text{diss}(u,\phi)$ vanishes when considering the analytical (smooth) solution (zero jumps), thus not altering the underlying physics in well resolved cases. We will investigate the effect of this non--linear dissipation through von Neumann analyses. These non--linearities arise as a result of the penalisation on the size of the interface jumps when increasing $\lambda$.

In the 3D Navier--Stokes equations \eqref{eq::NS1}, inviscid fluxes or Riemann solvers are usually constructed as the average of both adjacent states, plus an interface dissipation that depends on the two states $(\vec{u}_L,\vec{u}_R)$ \cite{2009:Toro}:
\begin{equation}
\vec{\boldsymbol{F}}_e^\star\cdot\vec{\boldsymbol{n}} = \{\!\!\{\vec{\boldsymbol{F}}\cdot\vec{\boldsymbol{n}}\}\!\!\} - \text{diss}(\vec{u}_L, \vec{u}_R).\label{eq:RiemannStructure}
\end{equation}
In this paper, we will consider Roe dissipation,

\begin{equation}
\text{diss}_{Roe}(\vec{u}_L, \vec{u}_R) = \sum_{e=1}^5 \alpha_e |\beta_e|\vec{K}_e,
\label{eq:RoeRiemann}
\end{equation}
where the intensities $\alpha_e$, eigenvalues $\beta_e$ and eigenvectors $\vec{K}_e$ are computed from the Roe averaged states \cite{2009:Toro}. We will study the effect of the parameter $\lambda$ in \eqref{eq:lambdastabdiss} in the Navier--Stokes Equations (NSE) by modifying \eqref{eq:RiemannStructure} to:

\begin{equation}
\vec{\boldsymbol{F}}^\star\cdot\vec{\boldsymbol{n}} = \{\!\!\{\vec{\boldsymbol{F}}\cdot\vec{\boldsymbol{n}}\}\!\!\} -\lambda \text{diss}(\vec{u}_L, \vec{u}_R),\label{eq:RiemannStructureWithLambda}
\end{equation}
to control the dissipation added through cell interfaces. This new expression \eqref{eq:RiemannStructureWithLambda} can now be compared to the linear advection form \eqref{eq:lambdastab}. This strategy has been already adopted to design low dissipation versions for Roe Riemann solvers by estimating an appropriate value for $\lambda$ \cite{2016:Obwald}.
Only Roe dissipation \eqref{eq:RoeRiemann} will be studied in this work.

\subsubsection{Second dissipation technique: discretisation of viscous terms} \label{subsubsec:theory_LES}

Numerical dissipation can also be added in the discretisation of the viscous terms.
Several techniques are available to discretise second order derivatives; a review can be found in \cite{2002:Arnold}.
Some of these methods, e.g. the Bassi--Rebay 1 (BR1) scheme, are neutrally stable \cite{2017:Gassner} (it adds the minimum dissipation required to achieve a stable scheme) whilst others may introduce some extra dissipation, e.g. the Symmetric Interior Penalty (SIP) scheme.
Only the BR1 scheme will be studied here, as it was proven in \cite{2018:Manzanero} that for 1D problems and Gauss--Lobatto nodes (also for 2D with cartesian meshes) the SIP scheme can be seen as a BR1 scheme plus a penalisation term with the same form as the introduced by upwind Riemann solvers. In fact, they are algebraically equivalent for an appropriate penalty parameter, $\sigma$, value.

Limited efforts have been devoted to understand the combined effect of upwind Riemann solvers and discretised viscous terms in previous works. Therefore, special attention will be paid in von Neumann analyses to the combined effect of the Riemann solver when applied to the advection--diffusion equation.
In the Navier--Stokes equations, we will compare the results obtained with von Neumann analyses to those with a constant physical viscosity since the former only allows the analysis for constant viscosities.

%

\subsubsection{Third dissipation technique: artificial viscosity (LES models and spectral vanishing viscosity)} \label{subsubsec:theory_SVV}

An alternative (or complementary) technique to add numerical viscosity to the numerical scheme is the one followed by artificial viscosity techniques. We will consider here two different approaches: a standard Smagorinsky LES model and spectral vanishing viscosity.

Briefly, LES models introduce a solution dependent viscosity, $\mu_{a}$, usually through a second order derivative. More specifically, the non--constant artificial viscosity $\mu_{a}(u)$ is estimated from the solution using sensors based on either numerical (e.g. solution jumps or modal sensors) \cite{2006:Persson,2010:Barter,2011:Klockner,2013:Persson} or physical models (e.g. LES models, or entropy and energy production) \cite{2014:Abbassi,2017:Kornelus}.

Precisely, we model $(f_{diss})_x$ in equation \eqref{eq:DGScheme} to introduce the artificial viscosity:
\begin{equation}
(f_{diss})_x = \frac{\partial}{\partial x}\biggl( \mu_{a} \frac{\partial u}{\partial x}\biggr).
\label{eq:advectionLES}
\end{equation}

Due to the linear nature of von Neumann analysis performed herein, we restrict ourselves to the effect of a constant viscosity $\mu_{a}$ and therefore the obtained results are equivalent to those obtained in the previous section (for the discretisation of viscous terms).

In the Navier--Stokes equations, the standard Smagorinsky subgrid--scale model \cite{1963:Smagorinsky}, is used in all studies shown in this paper. Precisely, $\vec{\boldsymbol{F}}_{diss}$ in \eqref{eq::NS5} is:
\begin{equation}
\vec{\boldsymbol{F}}_{diss} =
\vec{\boldsymbol{F}}_{v}(\mu=\mu_S,\vec{v},\nabla\vec{v}),
\end{equation}
where $F_{v}(\mu,\vec{v},\nabla\vec{v})$ is the Navier--Stokes viscous flux defined in \eqref{eq:viscousfluxes}. The Smagorinsky viscosity is defined as:
\begin{equation}
\mu_{S} = C_S^2\Delta^2 |\boldsymbol{S}|.
\end{equation}
The classical value for isotropic turbulence $C_S=0.2$ is selected whilst the filter width $\Delta$ is computed as in \cite{2017:Flad}:

\begin{equation}
\Delta^3 = \frac{\text{Cell volume}}{(N+1)^3},
\end{equation}
thus accounting for both the mesh elements size and the polynomial order.
The effect of the Smagorinsky model will be analysed in detail in Section \ref{sec:LESSVV_hybrid}.

%
%


%
%

The last approach to introduce numerical dissipation studied in this paper is the Spectral Vanishing Viscosity (SVV). This technique was initially introduced to stabilise  Fourier spectral methods in \cite{1989:Tadmor}, and later adapted to high--order continuous Galerkin methods in \cite{2000:Karamanos}. This technique considers a constant viscosity, $\mu_{SVV}$, which is applied unevenly on the different modes that form the solution. The operator that chooses the intensity of each mode in the dissipation is called the viscous kernel $Q_\mu$. Precisely, we model $(f_{diss})_x$ in equation \eqref{eq:DGScheme} to introduce the SVV regularisation:
\begin{equation}
(f_{diss})_x = \frac{\partial}{\partial x}\biggl( \mu_{SVV}Q_\mu \star \frac{\partial u}{\partial x}\biggr),
\label{eq:advectionSVV}
\end{equation}
where the operator $\star$ denotes the modal convolution operator, applied to the solution derivative and the SVV viscous kernel (see implementation details in \cite{2000:Karamanos}).

The spectral distribution of the viscosity is defined \textit{ad-hoc} in the viscous kernel, $Q_\mu$, of the SVV method. The most popular choice is the exponential distribution:

\begin{equation}
Q_\mu(k) = e^{(k-N)^2/(k-M)^2} \text{, if }k>M,~~ Q_\mu(k) = 0 \text{ otherwise,}
\end{equation}
where $M$ is a constant that sets the mode filter cut-off. Nevertheless, in this work we adopt the kernel introduced recently by Moura et. al \cite{2016:Moura}, who considered a potential law:

\begin{equation}
Q_\mu(k) = (k/N)^{P_{SVV}}, \quad k=0,1,\ldots,N. \label{eq:SVVKernel}
\end{equation}
In \eqref{eq:SVVKernel}, the constant $P_{SVV}$ is the kernel power coefficient. This kernel is convenient since a standard viscous discretisation is recovered when $P_{SVV}=0$, and because the dissipation vanishes if $P_{SVV}>>1$.
The SVV has been widely used and studied for continuous Galerkin and Fourier spectral methods \cite{2001:Tadmor,2006:Kirby,2016:Moura,2017:Ferrer}. In \cite{2016:Moura}, Moura et al. performed the dispersion--dissipation analysis of the continuous Galerkin SVV, finding similarities between the numerical dissipation introduced to that obtained using upwind Riemann solvers in discontinuous Galerkin methods. However, for discontinuous Galerkin methods there are not SVV studies available. One reason would be that the SVV was introduced to achieve similar dissipation behaviour in CG (where schemes present lack of numerical dissipation) to that obtained with upwind Riemann solvers in DG. Therefore, it may be intuitively thought that it is pointless to introduce an SVV method in DG, since the upwind fluxes achieve similar results, resulting in a more efficient and simpler implementation. However, in this paper we show that the SVV method is useful not only as a substitute of upwind Riemann solvers, but as a complement of the latter to control energy accumulations in medium wave--numbers.

In the Navier--Stokes equations, the SVV method is implemented by defining the following dissipation operator in \eqref{eq::NS5}:
\begin{equation}
\vec{\boldsymbol{F}}_{diss} = \left[\begin{array}{ccc}0 & 0 & 0\\
\hat{\tau}_{xx} & \hat{\tau}_{xy} & \hat{\tau}_{xz} \\
\hat{\tau}_{yx} & \hat{\tau}_{yy} & \hat{\tau}_{yz} \\
\hat{\tau}_{zx} & \hat{\tau}_{zy} & \hat{\tau}_{zz} \\
\sum_{j=1}^3 v_j\hat{\tau}_{1j} + k\hat{T}_x& \sum_{j=1}^3 v_j\hat{\tau}_{2j} + k\hat{T}_y& \sum_{j=1}^3 v_j\hat{\tau}_{3j} + k\hat{T}_z
\end{array}\right],
\label{eq:SVVNS}
\end{equation}
where hatted variables represent those with the SVV operator applied, e.g.:

\begin{equation}
\hat{\tau}_{ij} = \mu_{SVV}Q_\mu\star\bigl[\partial_jv_i + \partial_iv_j\bigr].
\label{eq::SVVNS_tau}
\end{equation}
%


\section{Numerical Results}\label{sec:VN_and_TGV}

We now study the three dissipation techniques introduced in the previous section by two means: performing a linear von Neumann analysis, and solving the Taylor--Green vortex problem. Von Neumann analysis gives a global vision of the behaviour of the linear advection--diffusion equation scheme in the wave--number domain independently of the initial condition. In contrast, the TGV problem, see \ref{subsec:theory_TGV}, allows to study the dissipation introduced numerically by inspecting the kinetic energy spectra of the flow. This section is organised with the same structure as Section \ref{subsec:DissipationTechniques}; we first study energy conserving schemes, to then analyse the three dissipation techniques introduced in Sections \ref{subsubsec:upwindRiemannTheory}, \ref{subsubsec:theory_LES} and \ref{subsubsec:theory_SVV}.

\subsection{Preliminaries on energy conserving schemes}\label{subsec:results_ECS}
We start by considering an energy conserving scheme. We study the constant speed advection equation with Gauss points and central fluxes (diffusion is not considered here). The choice of Gauss points for this von Neumann analysis is adopted since it represents a more traditional and widely used approach, but we found that all conclusions obtained hold with Gauss--Lobatto points (see \cite{2011:Gassner} for details).

\subsubsection{Von Neumann analysis}
\begin{figure}
	\centering
	\subfigure[Dispersion error.]{\includegraphics[scale=0.23]{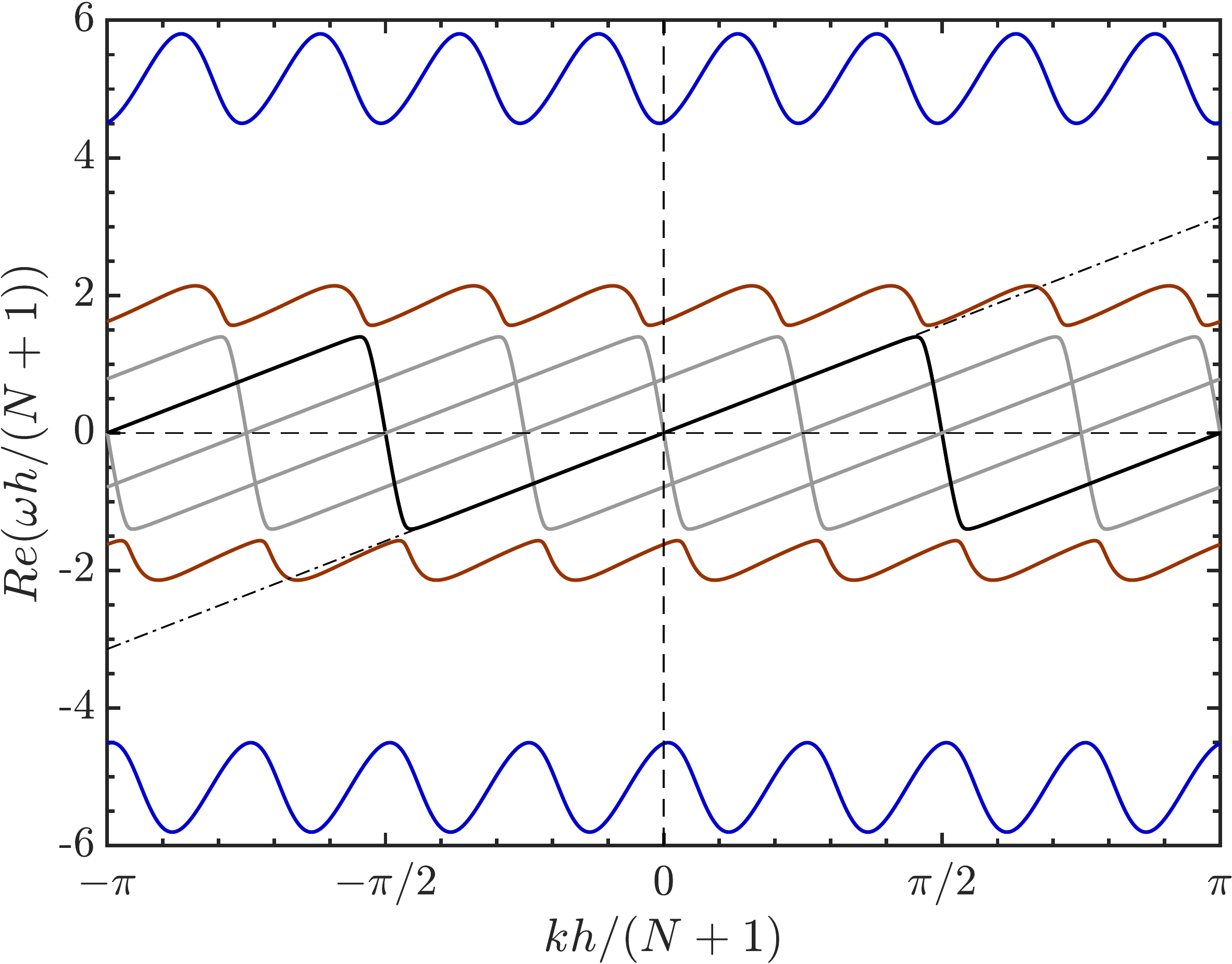}\label{fig:baseline_dispersion}}
	\subfigure[Diffusion error.]{\includegraphics[scale=0.23]{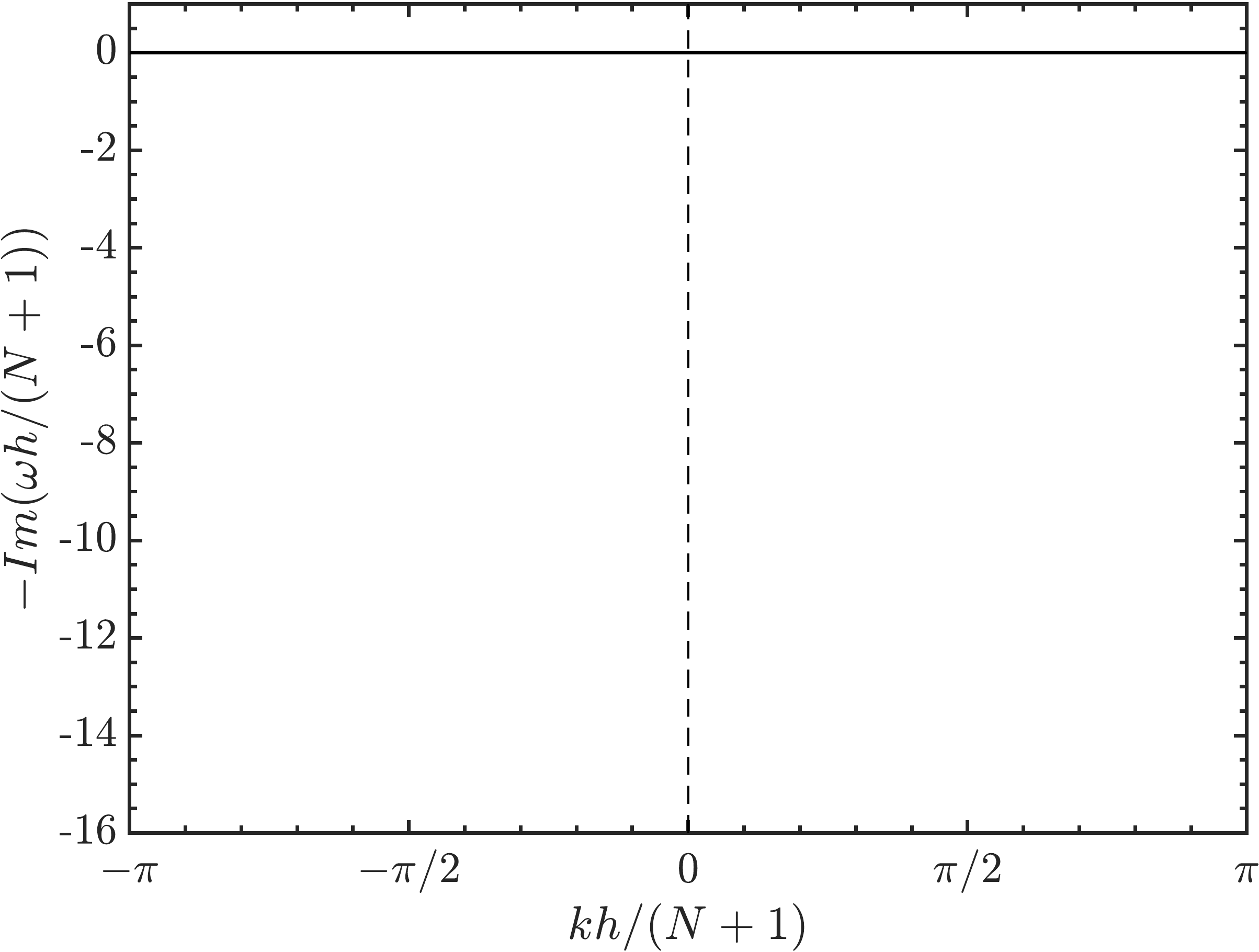}\label{fig:baseline_dissipation}}
	\caption{Dispersion-dissipation errors in the advection equation. Figures depict the eigenvalues obtained with an energy conserving DG (i.e. with central fluxes, $\lambda=0$) and polynomial order $N=7$. }
	\label{fig:baseline}
\end{figure}

Numerical dispersion and dissipation errors are depicted in Figures \ref{fig:baseline_dispersion} and \ref{fig:baseline_dissipation} respectively. The primary mode has been highlighted with a black line, whilst the rest are secondary modes. The dashed line in Figure \ref{fig:baseline_dispersion} depicts the analytical PDE speed (i.e. the straight line $Re(\omega)=k$). There are three groups of secondary modes according to their behaviours: first, grey modes in Figure \ref{fig:baseline_dispersion} are an exact traslation of the primary mode, therefore, they do not predict accurately the advection speed, and consequently this introduces numerical errors in the solution. Second, brown modes are medium-frequency modes that incur numerical errors in the propagation speed, except a narrow region close to $kh/(N+1)\simeq \pm \pi/2$, where these modes follow the analytical propagation speed, playing the role of the primary mode in that range. Third, blue lines represent non--physical high frequency modes. These results are in accordance with those previously reported in \cite{2008:Hesthaven}.

Figure \ref{fig:baseline_dissipation} shows dissipation errors, where neither numerical energy decay nor growth are experienced (i.e. the scheme is energy conserving, $Im(\omega)=0$). This result is consistent with numerical energy estimates: the scheme is energy conserving since the discrete energy balance inside the computational domain vanishes (when considering constant advection speeds, $a$). Further details can be found in \cite{Kopriva2} and also \cite{2017:Manzanero} for the generalisation to non--constant advection speeds.

\begin{figure}
\centering
\includegraphics[scale=0.3]{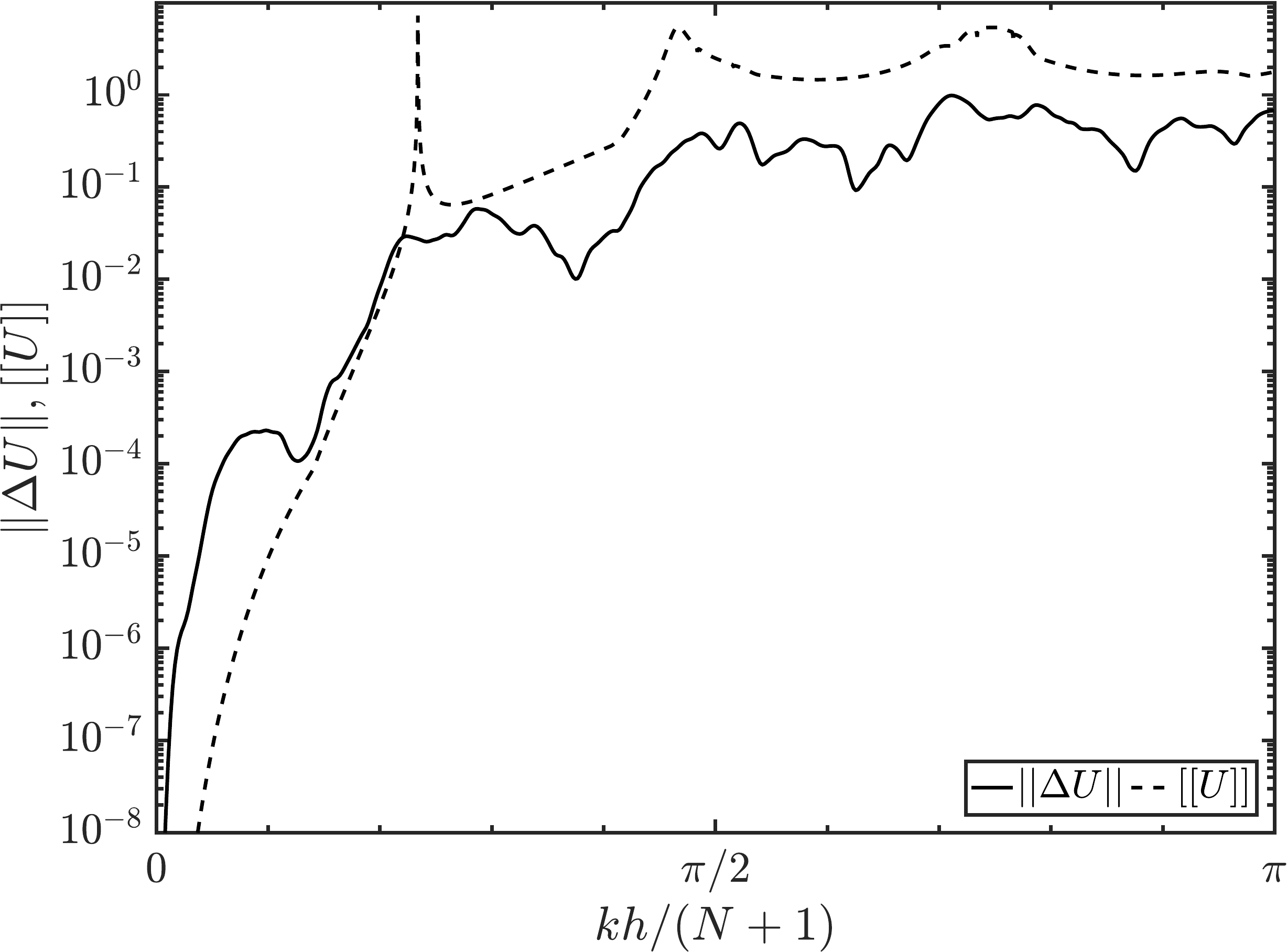}
\caption{Secondary modes error (solid line) of the energy conserving DG (Gauss points with central fluxes $\lambda=0$) and interface jumps (dashed line). We distinguish the two regions, first an exponential decay zone, where the solution is well-resolved, and an approximately flat region where both secondary mode errors and interface jumps reach their maximum values (where the solution is considered under--resolved). The interface jumps are represented alongside secondary modes error since in Section \ref{subsec:results_upwindRS} we find the former to be an accurate estimation of the latter.}
\label{fig:baseline_secondarymodes}
\end{figure}

Regarding the scheme accuracy, and following \eqref{eq:secondary_modes_error}, the influence of each mode depends on its projection onto the initial condition, $A_m$. To illustrate this, we have represented in Figure \ref{fig:baseline_secondarymodes} the secondary mode errors, $\Vert \Delta U\Vert $, against the non--dimensional wave--number, $\hat{k}$. This curve shows exponential decay for small wave--numbers, and an approximately flat region for high wave--numbers. This suggests a relationship between secondary mode errors and well/badly resolved regions.

\subsubsection{Navier--Stokes TGV problem}
\begin{figure}
	\centering
	\includegraphics[scale=0.4]{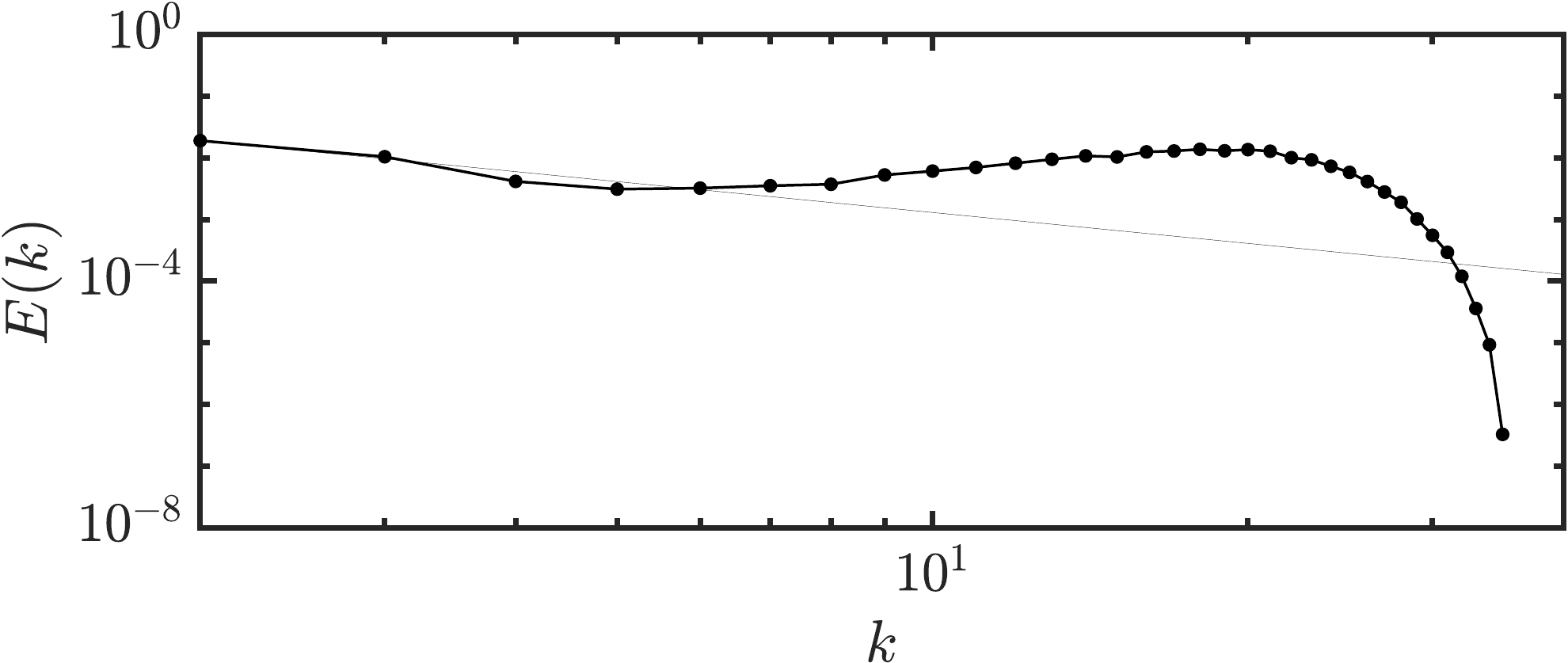}
	\caption{Kinetic energy spectra of the inviscid Taylor--Green vortex problem (with Mach number $M_0=0.1$) in $t=14$, alongside the theoretical Kolmogorov's solution $k^{-5/3}$ \cite{pope2001turbulent}. This result has been obtained with the energy conserving DG scheme introduced in \cite{2016:gassner}. We find that the lack of dissipation leads to energy accumulation in high wave--numbers, where the solution is severely under--resolved.}
	\label{fig:numExps:noDissipTGV}
\end{figure}

We show the inviscid Taylor--Green vortex solution (with Mach number $M_0=0.1$), see \ref{subsec:theory_TGV},  computed with a kinetic energy conserving scheme (detailed in \ref{appendix::A}. We consider a Cartesian $8^3$ mesh and fourth order ($N=4$) polynomials. Figure \ref{fig:numExps:noDissipTGV} depicts the kinetic energy spectra after $14$ time units. We find an undesired accumulation of energy in high wave--numbers as a result of the undissipated kinetic energy transferred from large to small eddies (i.e. the scheme is energy conserving). The solution is severely under--resolved, where high wave--number modes energy (with large dispersion errors) are not dissipated. Notice that this can be inferred with von Neumann analysis (see Figure \ref{fig:baseline}). A solution that presents large number of high wave--number spectral components, sees an accumulation of energy since, at these wave--numbers, dispersion errors are important and there is no dissipation to drain energy at these wave--numbers. Moreover, as a result of the insufficient energy drain, this accumulation is also transferred to low wave--numbers. Following sections study several techniques to introduce numerical dissipation and tackle this problem.

\subsection{First dissipation technique: upwind Riemann solver}\label{subsec:results_upwindRS}

In this section we study the effect on the scheme stability and accuracy of upwind Riemann solvers. This is performed by introducing the operator $\lambda\text{diss}(u,\phi)$ defined in \eqref{eq:lambdastab} in the PDE weak form \eqref{eq:DGScheme}. The motivation behind using a penalisation proportional to the solution interface jumps is that we can show that the amplitude of secondary modes is intimately related to the magnitude of solution interface jumps. To do so, we rewrite interface jumps (for instance, in the right boundary of an element and considering Gauss--Lobatto points for the sake of simplicity) as:

\begin{equation}
\llbracket u^{el}\rrbracket \bigr|_1 = u^{el}(1) - u^{el+1}(-1)= u^{el}(1) - e^{ikh}u^{el}(-1) =
u^{el}_N - e^{ikh}u^{el}_0.
\label{eq:jumps}
\end{equation}
Replacing \eqref{eq:secondary_modes_error} in \eqref{eq:jumps}:

\begin{equation}
\llbracket u^{el}\rrbracket \bigr|_1 = \boldsymbol{u}^{el}_{0,N}\exp[-i\omega_p t] + \Delta u_N -  \boldsymbol{u}^{el}_{0,0}e^{ikh}\exp[-i\omega_p t] - \Delta u_0 =\Delta u_N- \Delta u_0,
\end{equation}
where the initial condition contribution to interface jumps vanishes (since it is considered to be smooth). Thus, we conclude that interface jumps are intimately related to secondary mode errors, hence, providing a real measure of whether the simulation is resolved or under--resolved. Moreover, several resolution sensors using interface jumps were developed in the past, see \cite{2010:Barter}. In Figure \ref{fig:baseline_secondarymodes}, for completeness, we have represented both the secondary modes errors and the solution interface jumps, showing the parallelism that exists between both error measures. These results were obtained with von Neumann analysis of the energy conserving scheme.

\subsubsection{Von Neumann analysis}

\begin{figure}
	\centering
	\subfigure[Dispersion error.]{\includegraphics[scale=0.23]{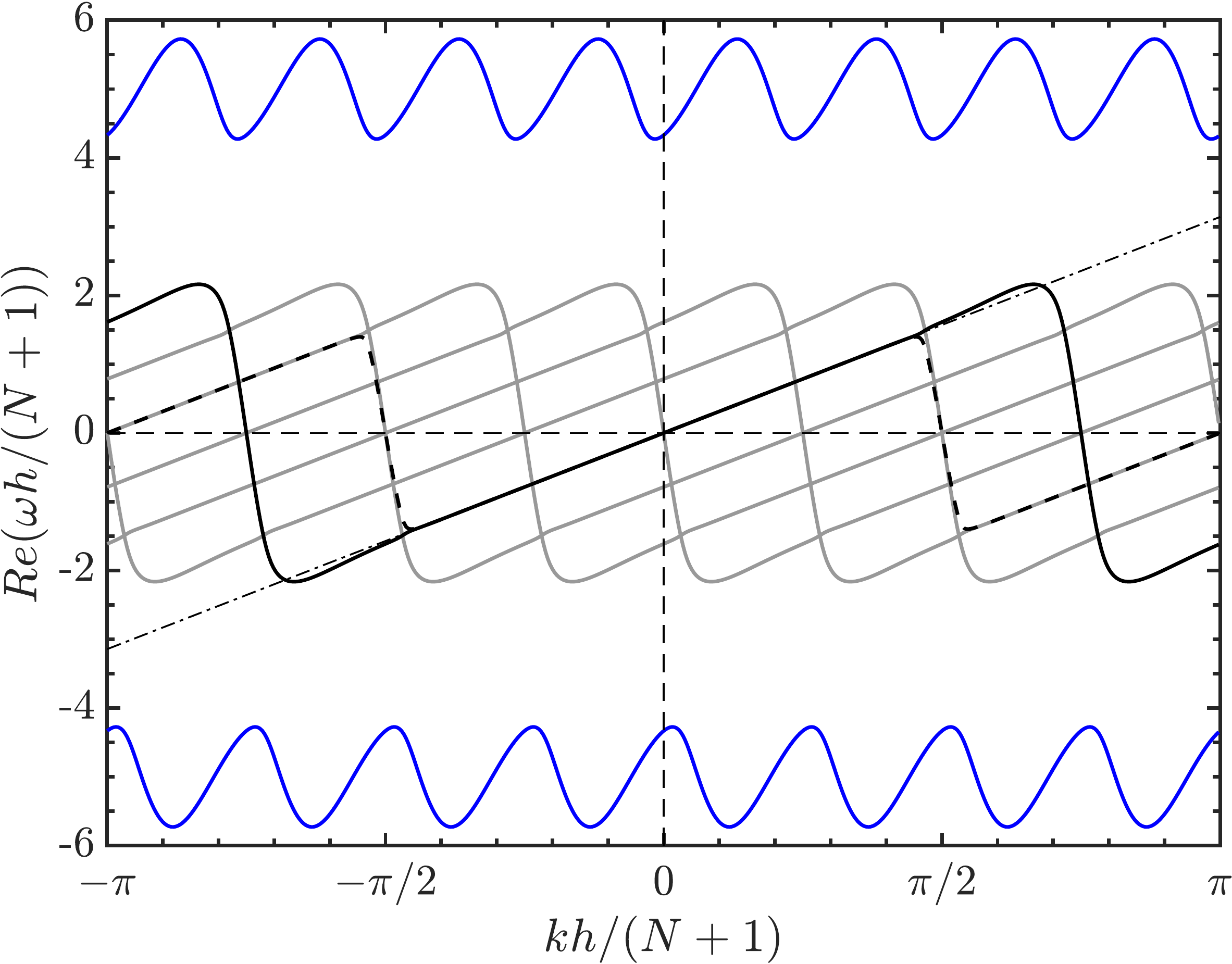}\label{fig:fluxes_0500_dispersion}}
	\subfigure[Dissipation error.]{\includegraphics[scale=0.23]{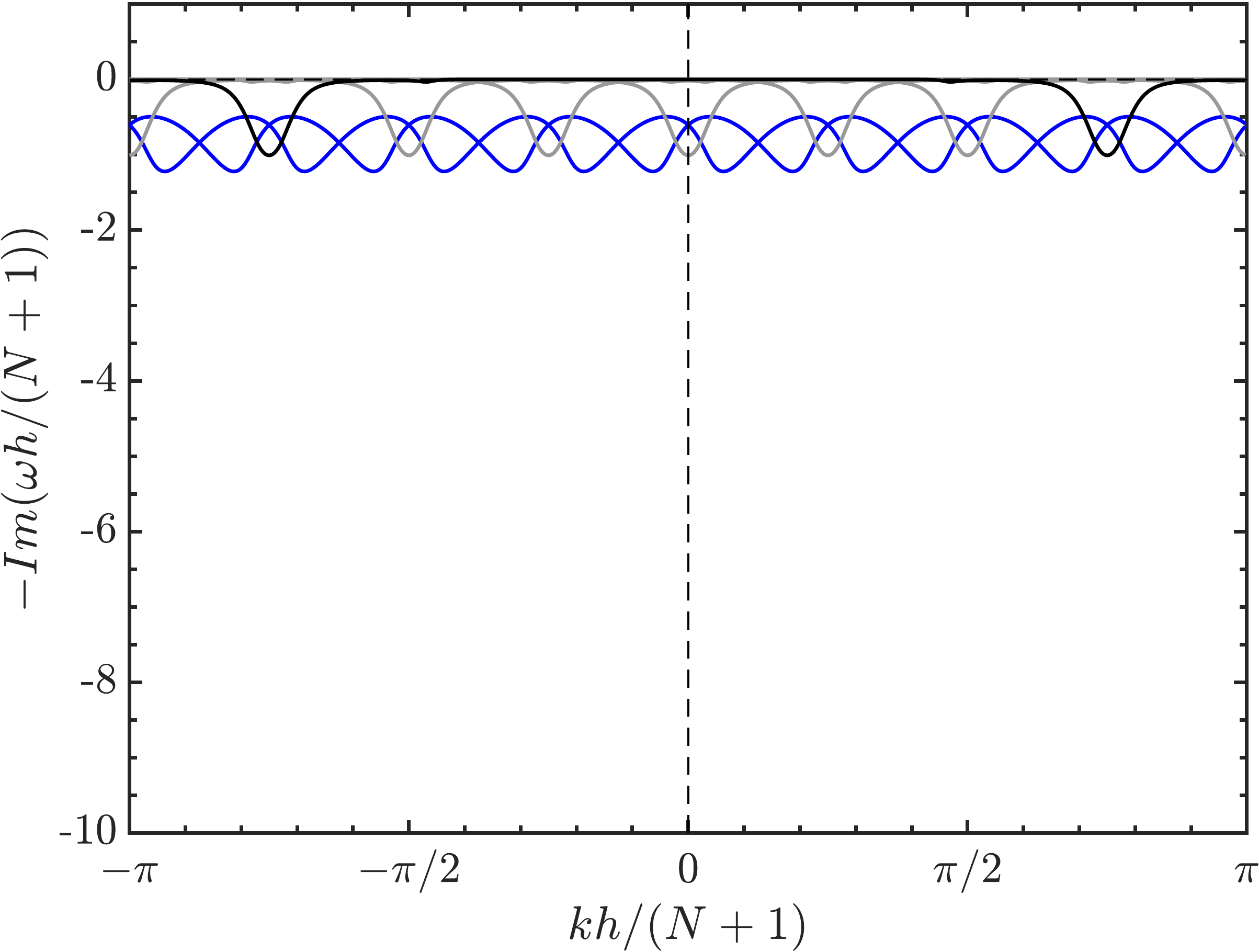}\label{fig:fluxes_0500_dissipation}}
	\caption{Dispersion and dissipation errors with $\lambda=0.5$.}
	\label{fig:fluxes_0500}
\end{figure}

\begin{figure}
	\centering
	\subfigure[Dispersion error.]{\includegraphics[scale=0.23]{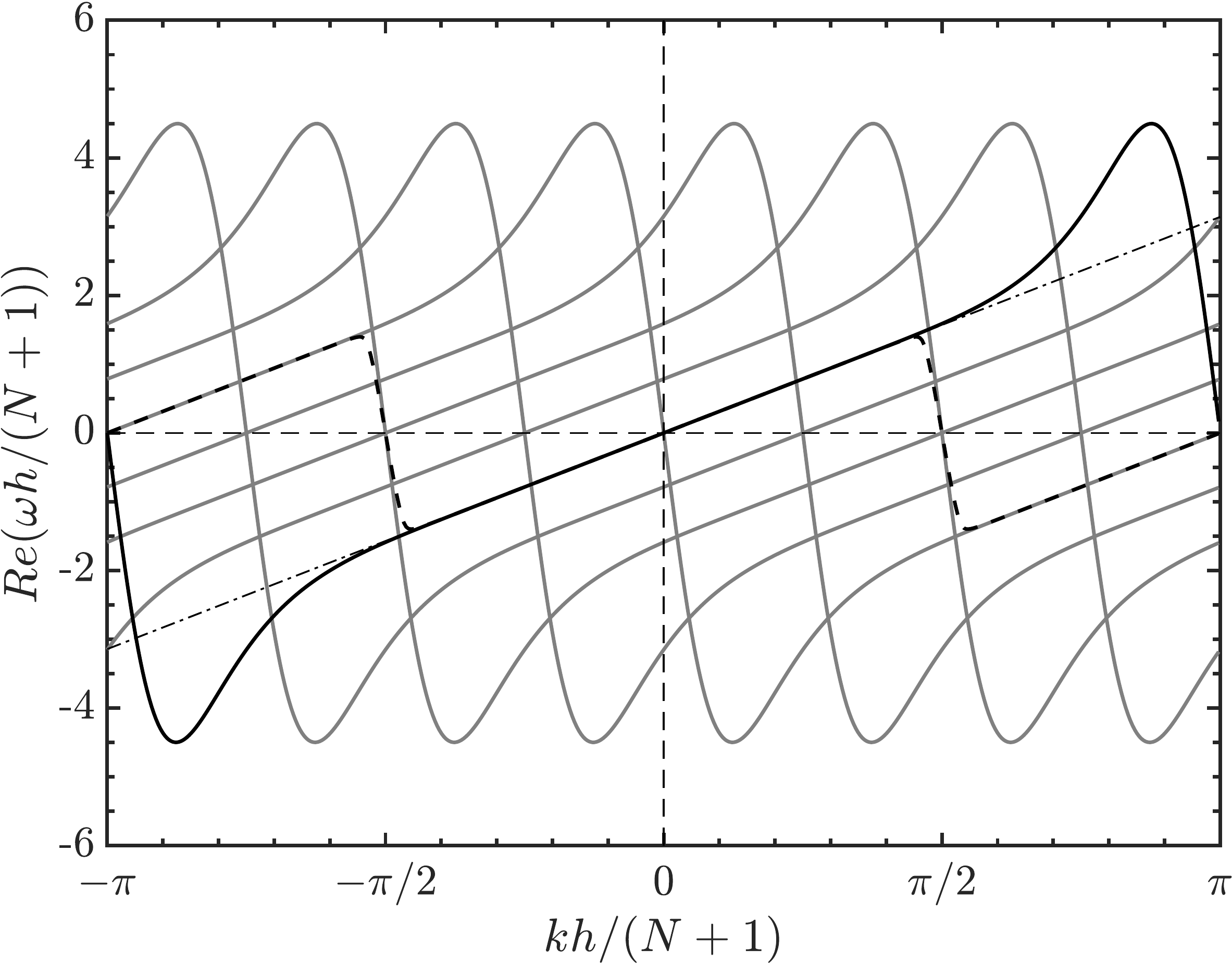}}
	\subfigure[Dissipation error.]{\includegraphics[scale=0.23]{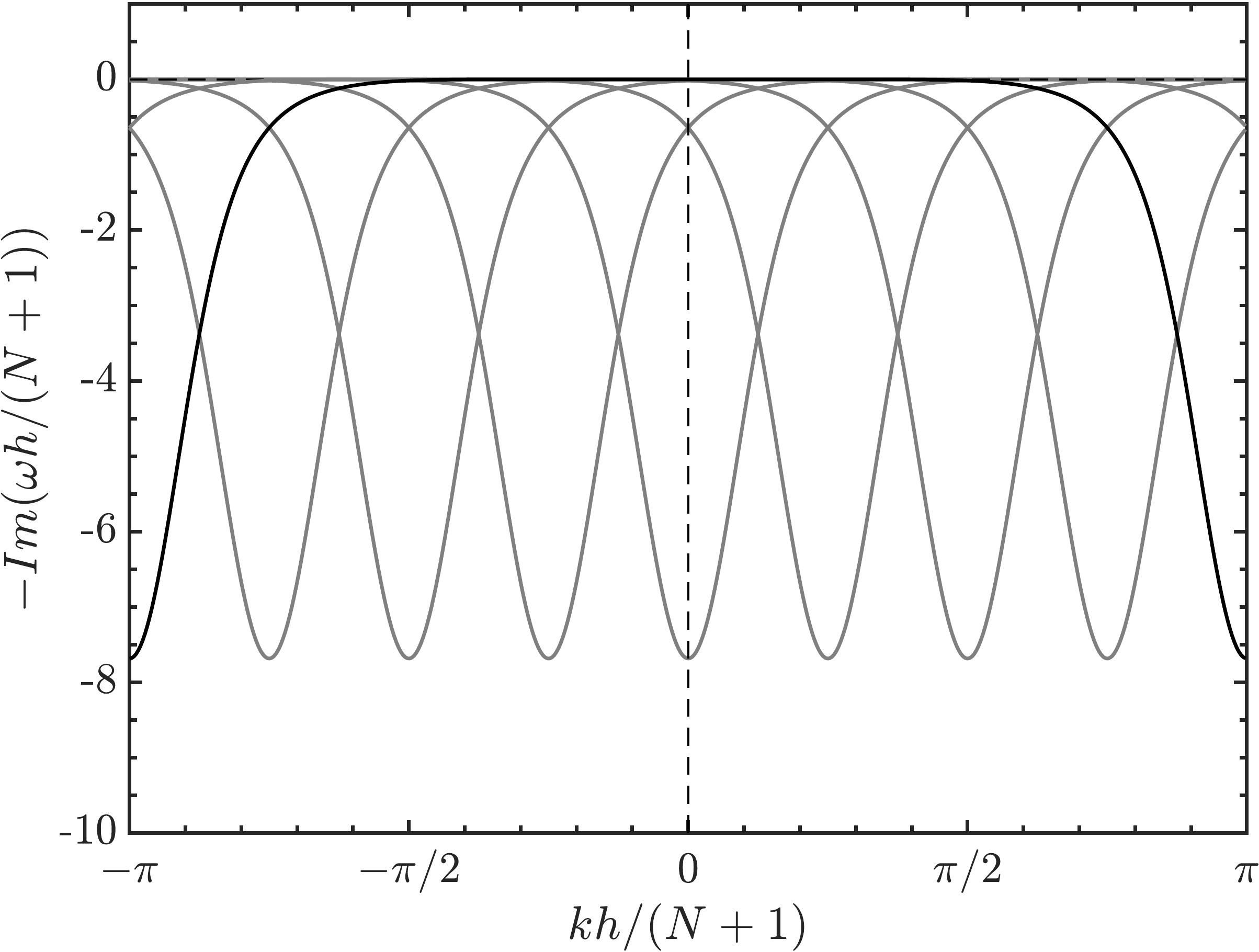}\label{fig:fluxes_1000_dissipation}}	
	\caption{Dispersion and dissipation errors with $\lambda=1.0$.}
	\label{fig:fluxes_1000}	
\end{figure}

\begin{figure}
	\centering
	\subfigure[Dispersion error.]{\includegraphics[scale=0.23]{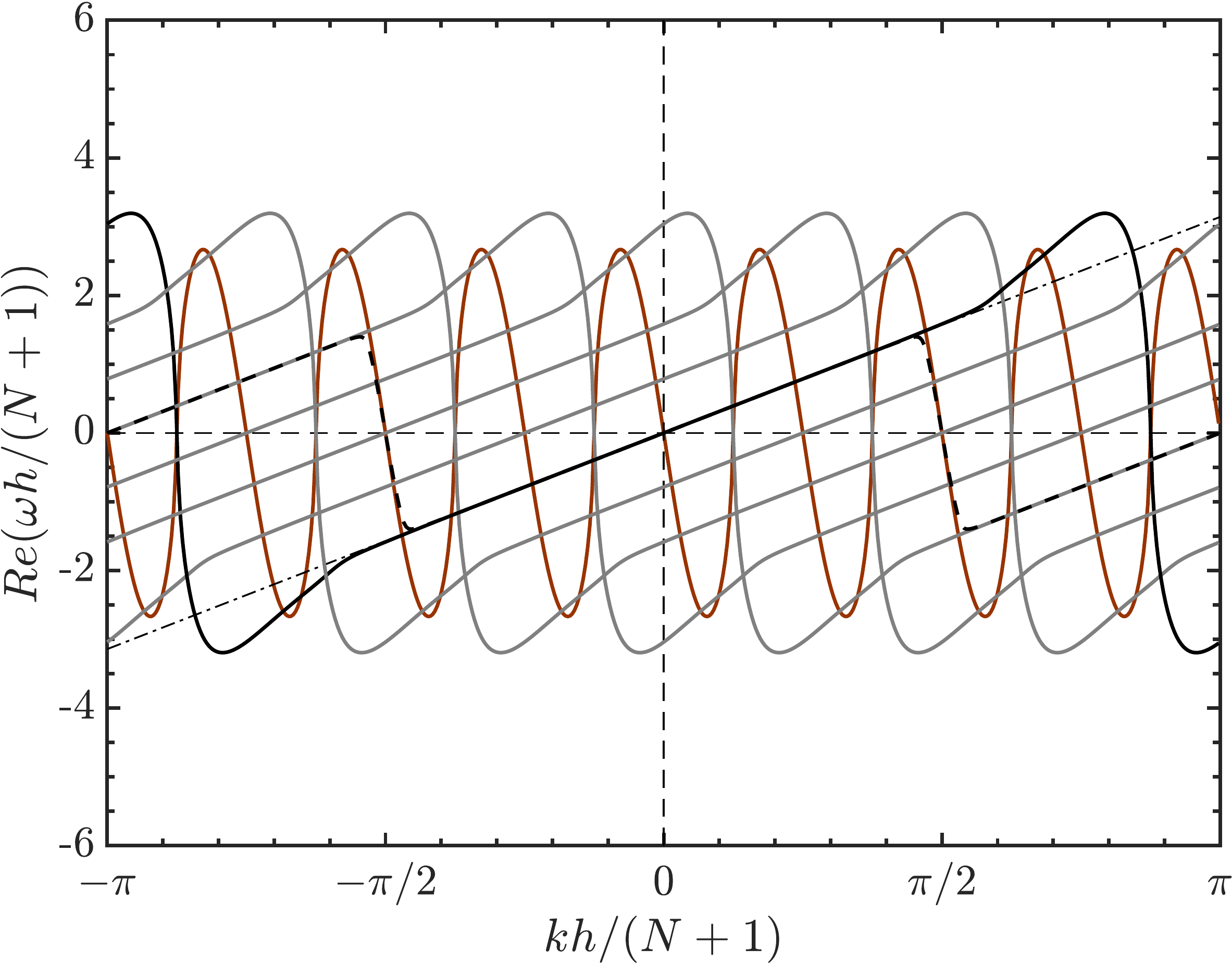}}
	\subfigure[Dissipation error.]{\includegraphics[scale=0.23]{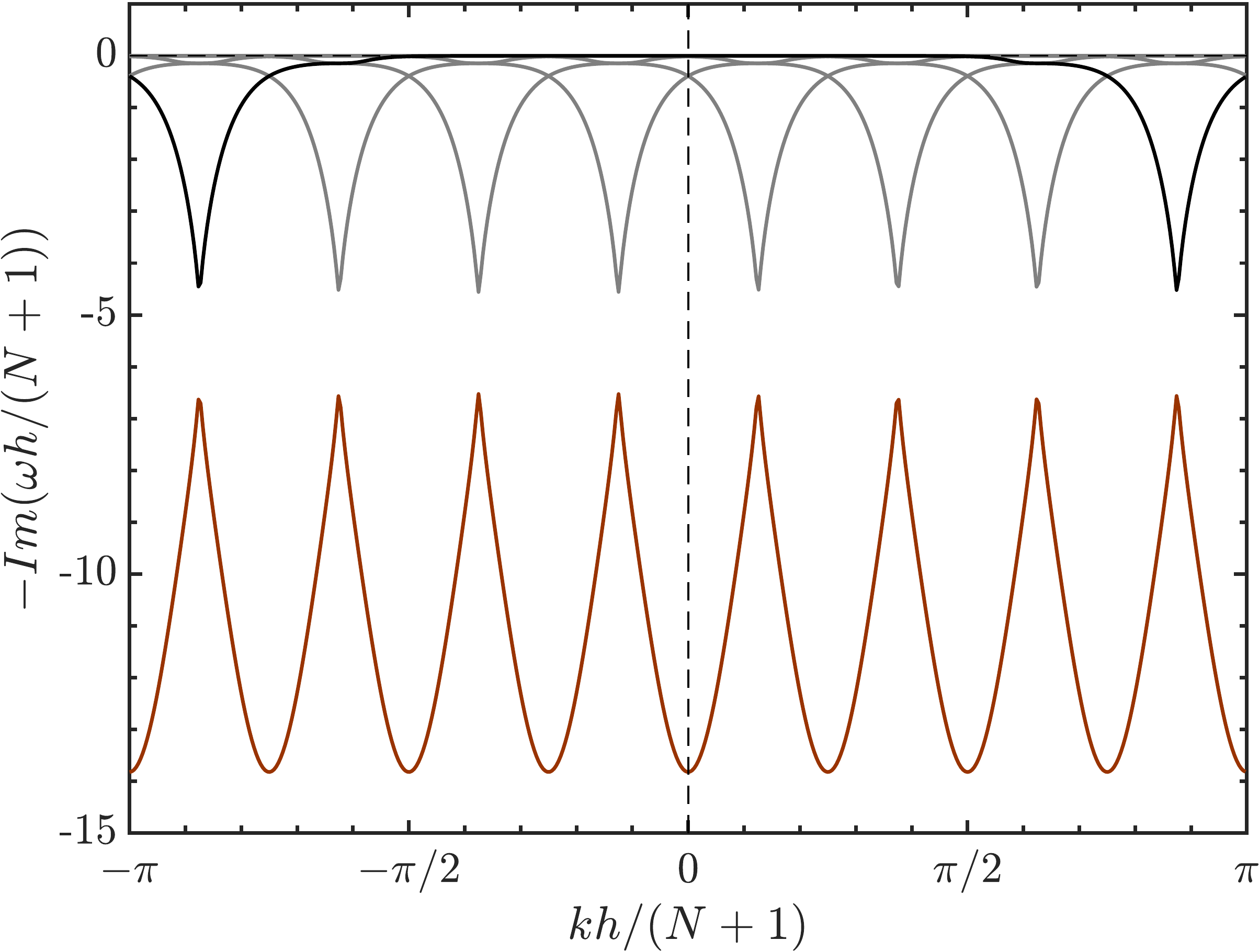}\label{fig:fluxes_1625_dissipation}}
	\caption{Dispersion and dissipation errors with $\lambda=1.625$.}	
	\label{fig:fluxes_1175}	
\end{figure}

 Dispersion and dissipation errors are depicted in Figures \ref{fig:fluxes_0500}-\ref{fig:fluxes_1175} for increasing values of $\lambda$.

 First, we show in Figure \ref{fig:fluxes_0500} the results with $\lambda=0.5$. We start by analysing dispersion errors (in Figure \ref{fig:fluxes_0500_dispersion}), where we distinguish two behaviours: the primary mode and its replications (black and grey lines), and two high frequency secondary modes, represented with blue lines. Compared to the  eigenvalues with $\lambda=0$ (in Figure \ref{fig:baseline_dissipation}), the brown medium-frequency modes have merged with the primary mode and its replications. For a direct comparison, the $\lambda=0$ primary mode is represented in black dashed line, whilst the solid black line depicts the $\lambda=0.5$ primary mode. When $\lambda=0.5$, the frequency range in which the primary mode remains accurate is extended, as a result of the merge with the brown modes.

 Regarding dissipation errors, we identify the same groups: the primary mode and its replications and the dissipation experienced by the high frequency modes (blue lines). Besides, the dissipation of the primary mode follows the traditional behaviour of high--order numerical dissipation: a selective filter, which introduces dissipation only at small scales (high wave--numbers). We conclude that high polynomial orders have an associated high bandwidth, as described in \cite{2008:Hesthaven,2015:Moura}.

When $\lambda=1.0$ we recognise a particular scenario, represented in Figure \ref{fig:fluxes_1000} (note that this value recovers upwind fluxes, reported in \cite{2015:Moura}. In this case, all modes follow the same behaviour, which might be regarded as optimal from the point of view of time-stepping limitations in explicit or   implicit-iterative solvers.

The results discussed so far indicate that the high frequency modes found with $\lambda=0$ and $\lambda=0.5$ have fully merged with the primary mode and its replications, increasing even more the accuracy of the latter to even higher wave--number ranges. This can be directly compared since the $\lambda=0$ primary mode has been also represented in black dashed line. The dissipation error, shown in Figure \ref{fig:fluxes_1000_dissipation}, shows the typical high--order filtering (similarly to $\lambda=0.5$), but being now all secondary modes (grey lines) replications of the primary mode (black line).

Lastly, when increasing $\lambda$ above $1.0$, the scheme changes again its behaviour. This has been represented in Figure \ref{fig:fluxes_1175}, for the value $\lambda=1.625$. Precisely, one mode (brown line) separates from the set, whose damping increases as $\lambda$ increases. On the contrary, the rest of the modes' (i.e. the primary mode and its replications) dissipation decrease, approaching zero as $\lambda$ increases. Thus, the eigenvalues tend to that of the continuous Galerkin method (whose dispersion--dissipation errors were studied in \cite{2016:Moura}) as $\lambda$ tends to infinity, with an extra mode (in brown line) whose numerical dissipation tends to infinity as $\lambda$ increases. This is consistent with the fact that the upwinding stabilisation $\lambda \text{diss}(u,\phi)$ acts as a penalty to interface jumps, and therefore is expected that the jumps tend to zero (i.e. continuous Galerkin approach) as $\lambda$ increases. However, this is achieved at the expense of an  infinitely large eigenvalue, and therefore, of a poor conditioning of the scheme.

\begin{figure}
	\centering
	\includegraphics[scale=0.45]{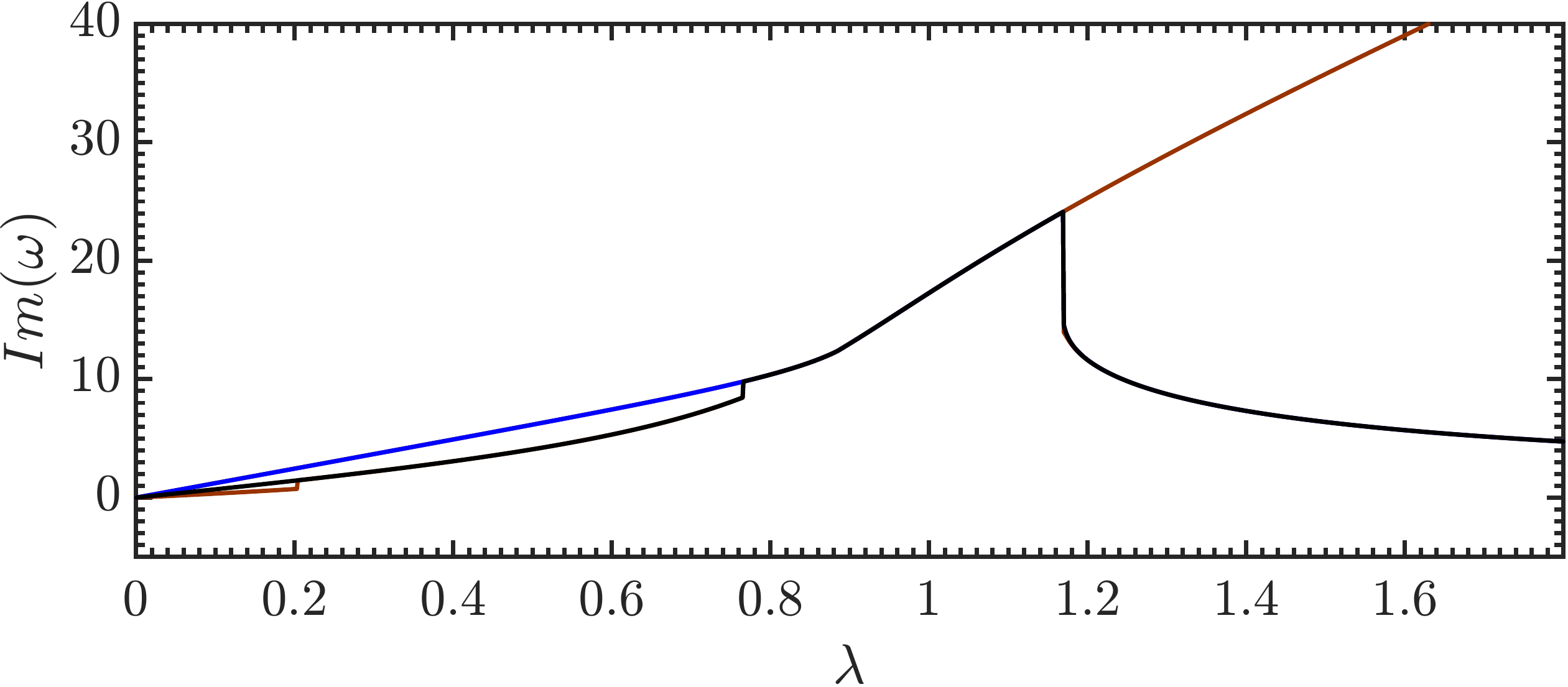}
	\caption{Diffusion of the different mode sets, and the effect of the Riemann solver parameter $\lambda$. This representation allows to follow the different bifurcation points that change the behaviour of the scheme dissipation. The colours represent the different secondary modes groups identified in Figures \ref{fig:fluxes_0500_dissipation}, \ref{fig:fluxes_1000_dissipation}, and \ref{fig:fluxes_1625_dissipation}. The precise values for the cases studied in the TGV (next section) are 0.7038 ($\lambda=0.1$), 17.27 ($\lambda=1$), 0.4065 ($\lambda=10$).}
	\label{fig:dissipation_lambda}
\end{figure}

Figure \ref{fig:dissipation_lambda} summarises the effect of upwind Riemann solvers $\lambda \text{diss}(u,\phi)$ in the numerical dissipation, where we depict the modes maximum dissipation rate for increasing $\lambda$ values. We start by considering central fluxes ($\lambda=0$), with three different set of dissipation-free modes, to observe that the scheme evolves to a state with the two set of modes in $\lambda\simeq 0.32$. Then, both collapse when $\lambda \simeq 0.77$. The effect of $\lambda$ in this region can be considered linear. Lastly, when $\lambda \simeq 1.17$, the set of modes separates into two: the primary mode and replication sets, and the dissipated mode that forces the solution continuity (i.e. approaches to a continuous Galerkin method behaviour).

\begin{figure}
	\centering
	\includegraphics[scale=0.3]{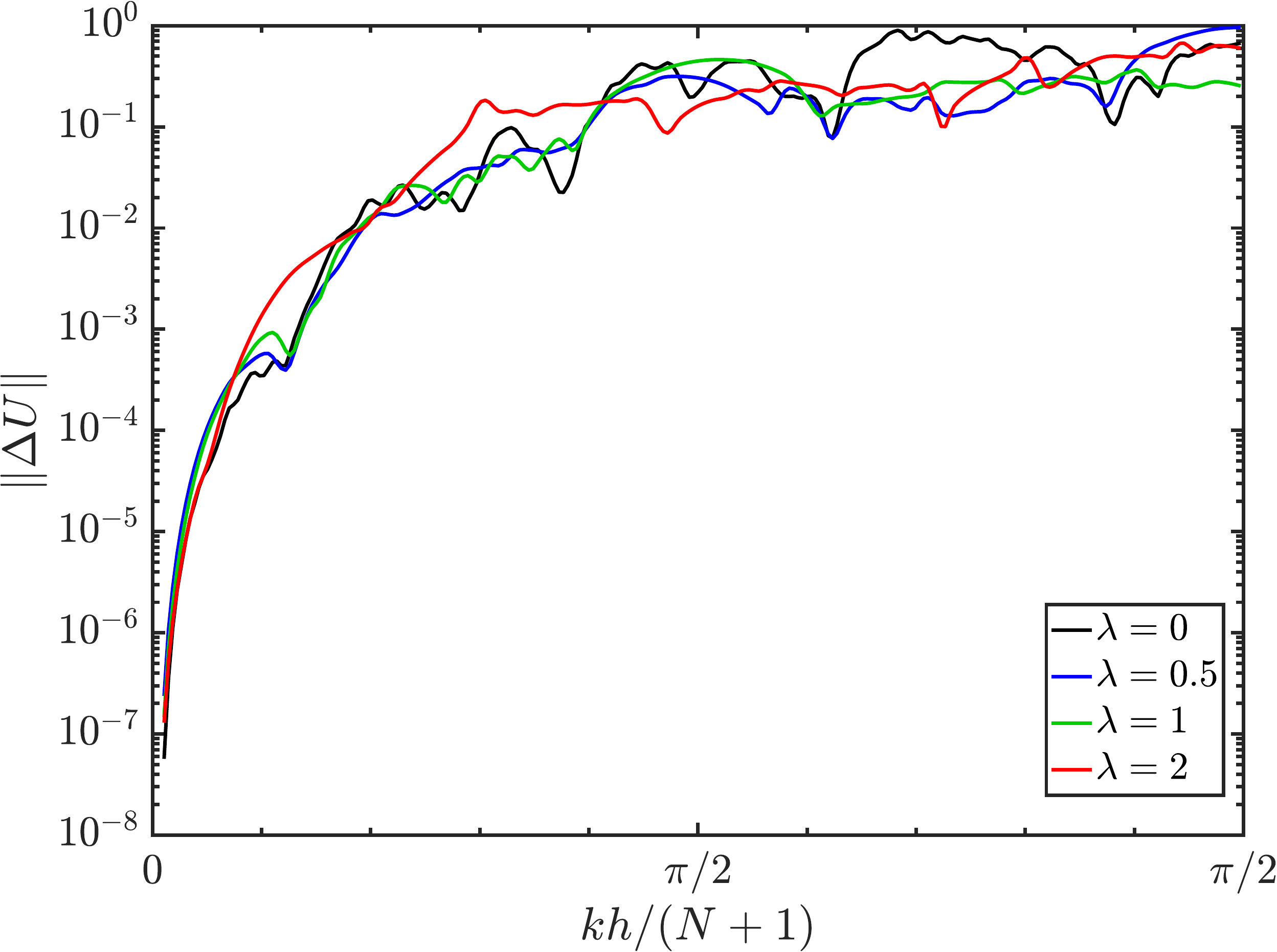}
	\caption{Upwind Riemann solver parameter $\lambda$ effect on secondary modes error.}
	\label{fig:fluxes_secondarymodes}
\end{figure}

To provide insight into the scheme accuracy, the secondary mode errors are represented against the dimensionless wave--number, for several $\lambda$ values in Figure \ref{fig:fluxes_secondarymodes}. It is found that high $\lambda$ values (i.e. above $1.0$) tend to increase the secondary mode errors in the low wave--number range. However, the effect of $\lambda$ does not seem to have a remarkable impact in the secondary modes contribution to the solution.

To summarise, upwind Riemann solvers do not have a significant effect on the solution accuracy, but they introduce dissipation (for moderate $\lambda\sim1$ values) where numerical errors (both dispersion and secondary mode errors) are important.

\subsubsection{Navier--Stokes TGV problem}
In the linear von Neumann analysis we have found that upwind Riemann solvers introduce dissipation in high wave--numbers. We have found that the amount of dissipation is non--linear with $\lambda$, and that values above approximately $\lambda=1$ yield lower dissipation rates. In this section, we study how this dissipation is introduced in the non--linear Euler equations. To do so, we solve the inviscid TGV problem (with Mach number $M_0=0.1$) introducing the $\lambda$ stabilisation based on Roe Riemann solver described in \eqref{eq:RoeRiemann}.

\begin{figure}
\centering
\subfigure[Numerical viscosity]{\includegraphics[scale=0.23]{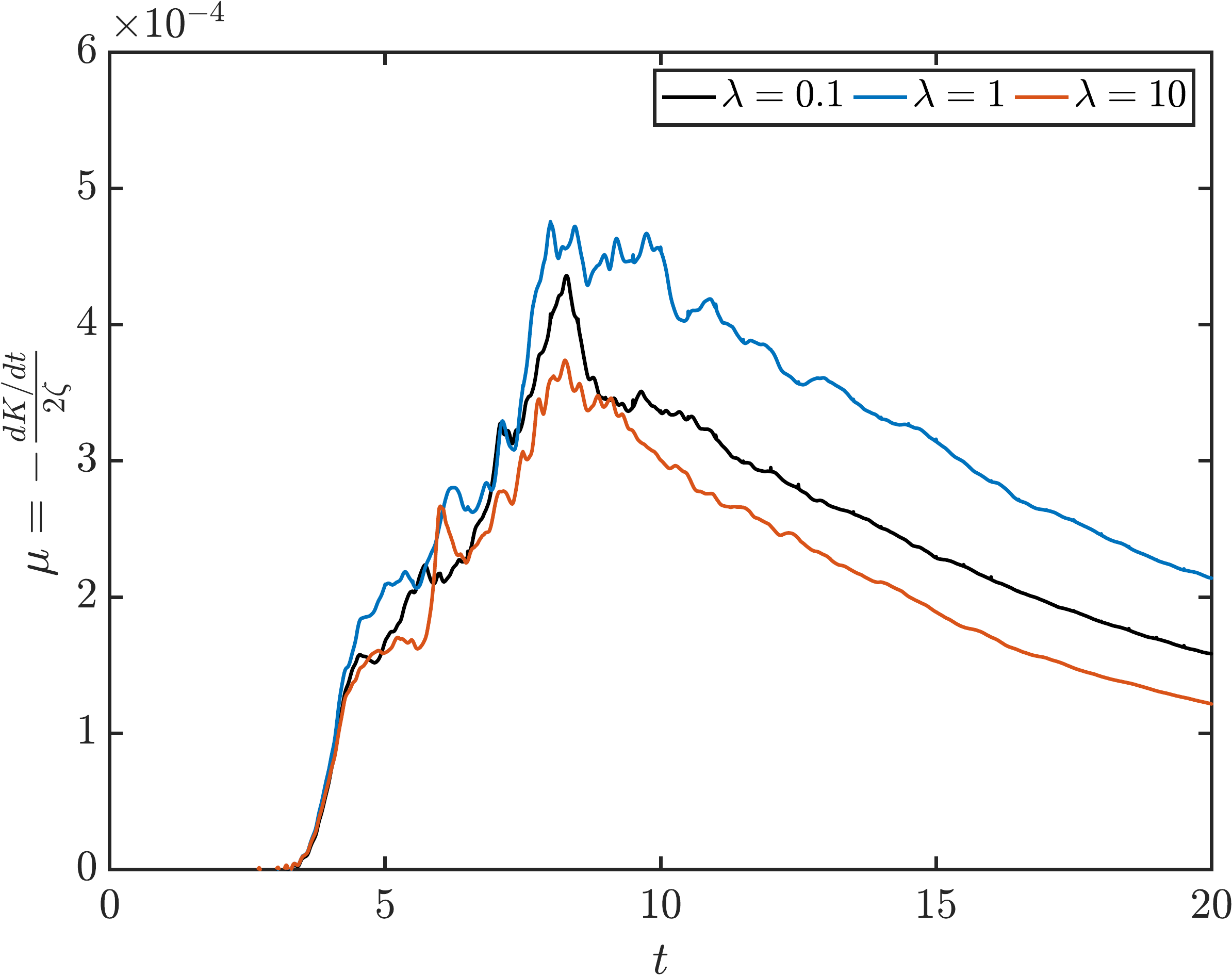}\label{fig:numExp:lambdaStudy:mu}}
\subfigure[Kinetic energy spectra in $t=14$]{\includegraphics[scale=0.23]{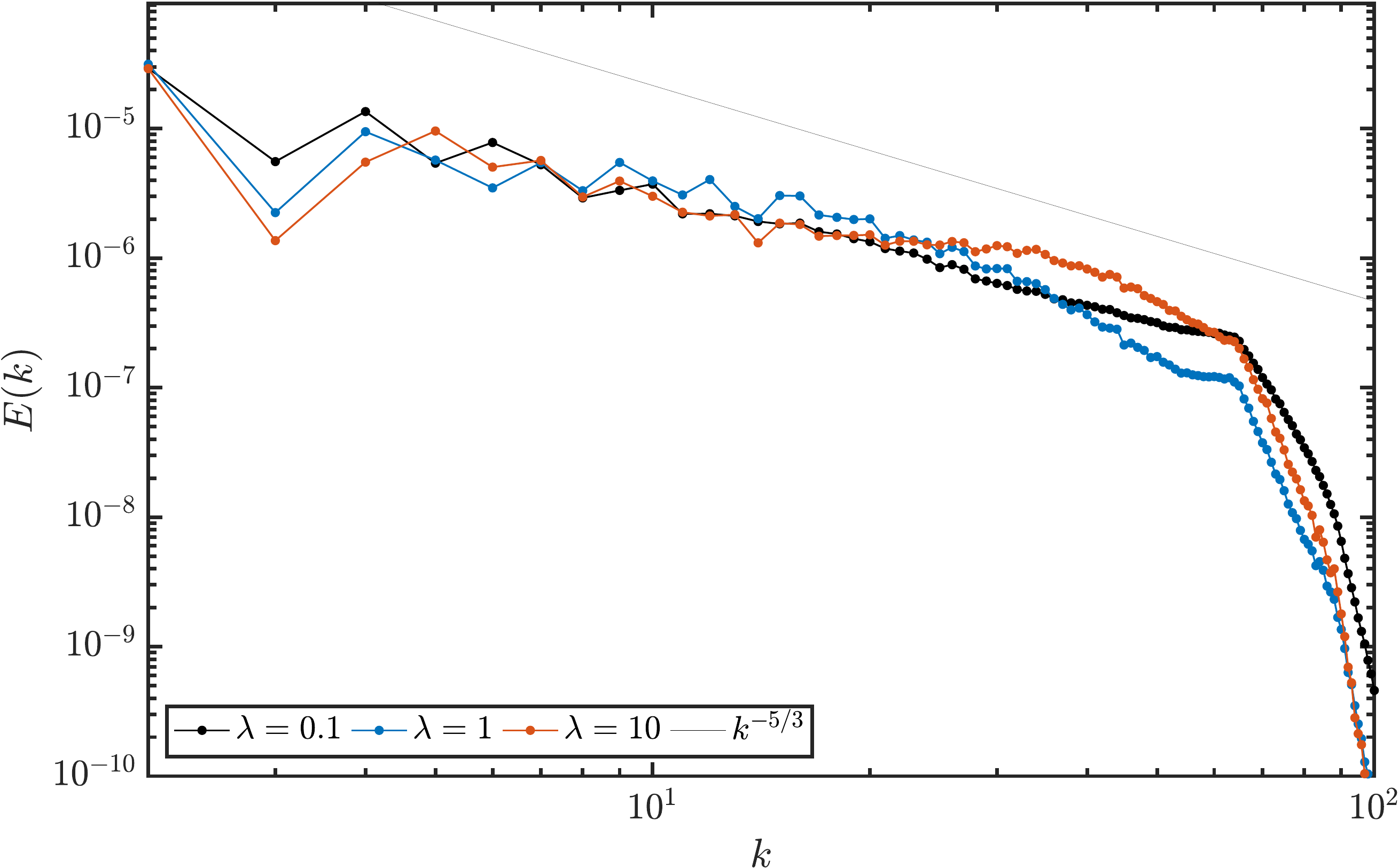}\label{fig:numExp:lambdaStudy:spectra}}
\caption{Inviscid Taylor--Green vortex problem ($M_{0}=0.1$). The configuration is a periodic box with $32^3$ elements and polynomial order $N=3$. For the interfaces, we have used Roe Riemann solver with the lambda stabilisation $\lambda \text{diss}(u,\phi)$ as defined in \eqref{eq:RiemannStructureWithLambda}. We have considered the values $\lambda=0.1,1,$ and $10$, whose results show clear parallelism with von Neumann analysis performed in Figure \ref{fig:dissipation_lambda}.}
\end{figure}
The structured mesh is constructed with $32^3$ elements, and the solution is approximated by $N=3$ polynomials.
Numerical dissipation for $t<20$ is depicted in Figure \ref{fig:numExp:lambdaStudy:mu}. We have considered three $\lambda$ values: 0.1 (low dissipation Roe), 1 (standard Roe), and 10 (hyper-upwind). Recall that the dissipation is only introduced numerically, since we consider the inviscid Euler equations. We find that the maximum dissipation is achieved with $\lambda=1$, which is in agreement with our von Neumann results in Figure \ref{fig:dissipation_lambda}. Therefore, for $\lambda\lesssim 1$ increasing $\lambda$ increases the scheme numerical dissipation, but increasing $\lambda$ for $\lambda\gtrsim 1$ yields the reverse effect, as demonstrated in von Neumann analysis (see Figure \ref{fig:dissipation_lambda}). In Figure \ref{fig:numExp:lambdaStudy:spectra} we show the kinetic energy spectra in $t=14$. We find that all three $\lambda$ values behave similarly in the low and medium wave--numbers range. At high wave--numbers, we find that the maximum dissipation is achieved by the standard Roe Riemann solver ($\lambda=1$), whose energy spectra shows an \textit{energy bump} (i.e. energy accumulation at high wave--numbers), originally reported in \cite{2018:Winters}. This configuration shows higher dissipation compared to the simulation with $\lambda=0.1$, and can be regarded an over-dissipated solution. Note that $\lambda=0.1$ is precisely the value provided by the low dissipation Roe Riemann solver derived in \cite{2016:Obwald} for $M_0=0.1$. Finally, Roe Riemann solver with $\lambda=10$ suffers  an accumulation of energy at high wave--numbers, as a result of its lack of dissipation compared to lower $\lambda$ values (see Figure \ref{fig:numExp:lambdaStudy:mu}).

\subsection{Second dissipation technique: discretisation of viscous terms}\label{subsec:results_LES}

In this section we study the effect of the discretisation of viscous terms. We start by discretising the advective part using central fluxes ($\lambda=0$, which entails zero dissipation)  such that we only focus on the effect of the discretisation of viscous terms. Subsequently, we discuss the effect of combined upwind Riemann solvers and viscous terms in the advection--diffusion equation.

\subsubsection{Von Neumann analysis: discretisation of viscous terms (viscosity $\mu$) and central fluxes ($\lambda=0$)}
This section is devoted to the assessment of the discretisation of viscous terms without adding extra stabilisation through upwind Riemann solvers. To do so, we perform von Neumann analysis of equation \eqref{eq:AdvEq} with Peclet number:

\begin{equation}
Pe=\frac{aL}{\mu}=1,
\end{equation}
where $L$ is the domain length. The dissipation introduced through this term is linear (and proportional to the inverse of $Pe$) in the advection--diffusion equation, according to the general analytical solution of \eqref{eq:AdvEq},

\begin{equation}
u(x,t) = e^{-k^2t/Pe}e^{i(kx-\omega t)}.\label{eq:AnalyticalSolutionWithPeclet}
\end{equation}

\begin{figure}
	\centering
\subfigure[Dispersion error.]{\includegraphics[scale=0.23]{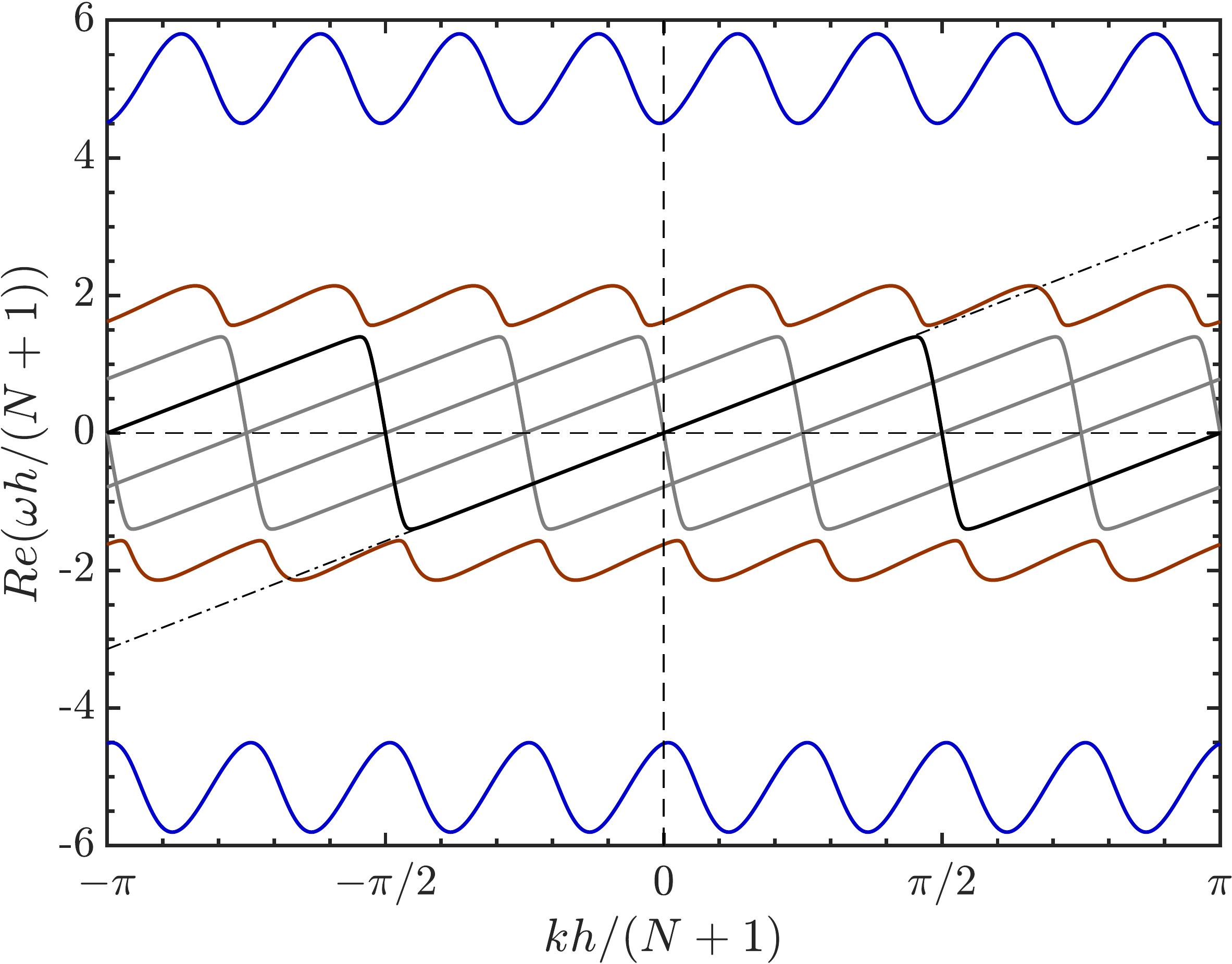}\label{fig:BR1_dispersion}}
\subfigure[Dissipation error.]{\includegraphics[scale=0.23]{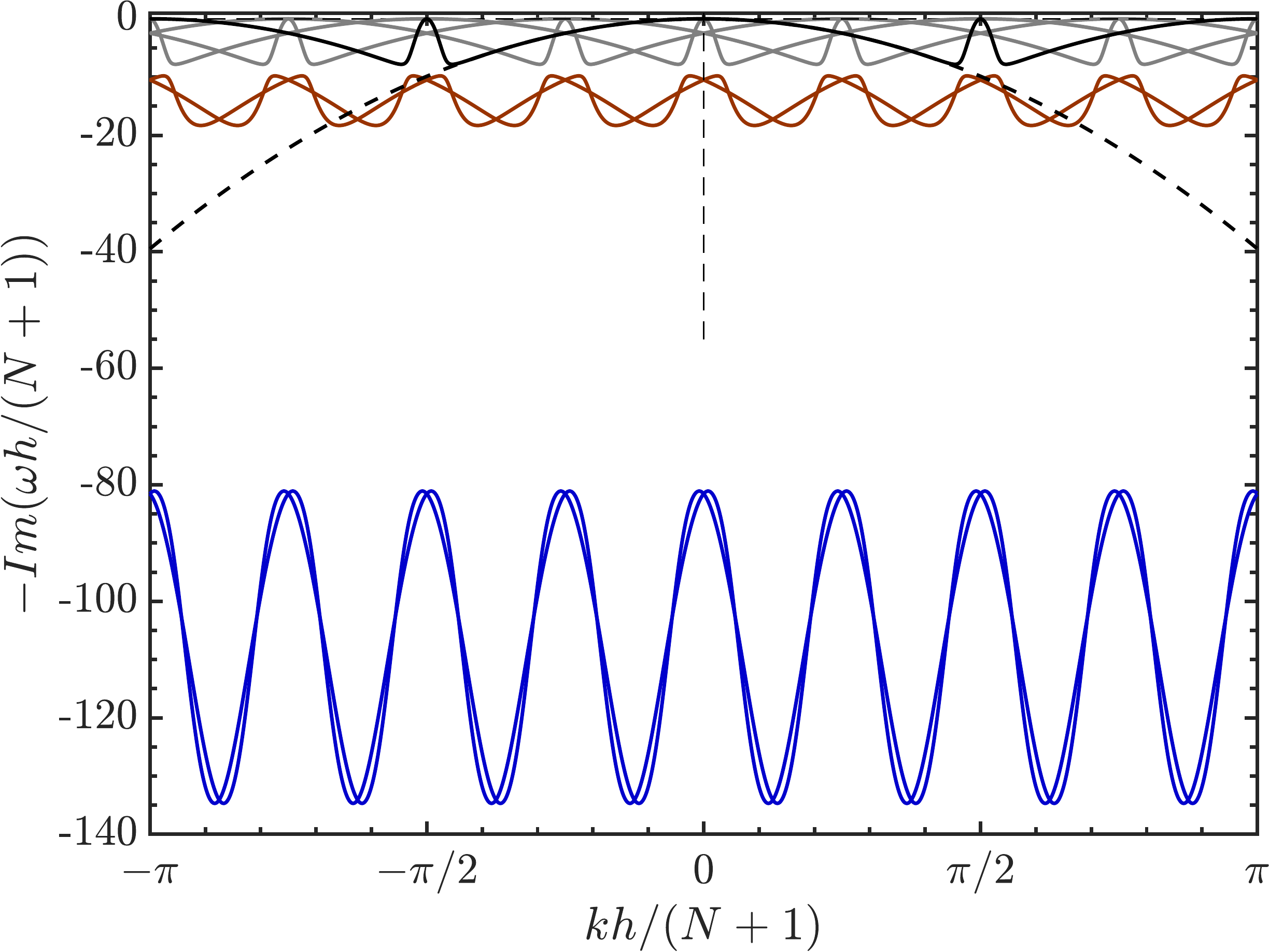}\label{fig:BR1_dissipation}}
\caption{Dispersion and dissipation errors obtained with the artificial viscosity (BR1 scheme) and central fluxes ($\lambda=0$). In this case, the Peclet number is $\text{Pe}=1$.}	
\label{fig:BR1_dispdiss}	
\end{figure}

In Figure \ref{fig:BR1_dispersion} we depict dispersion errors. We find that when adding a viscous discretisation with the BR1 formulation, dispersion errors remain unchanged with respect to the energy conserving scheme (Figure \ref{fig:baseline_dispersion}). In other words, dispersion errors generated by advection terms discretisation, and dissipation errors as a result of the discretisation of viscous terms are decoupled, and thus, the numerical dissipation introduced also varies linearly with the inverse of $Pe$.
As in the energy conserving scheme, we distinguish the primary mode and its replications (black and gray lines), medium-frequency modes (brown lines) and high frequency modes (blue lines). The dissipation carried by these modes, once the  viscosity is introduced, is represented in Figure \ref{fig:BR1_dissipation}. High-frequency modes are highly damped, and medium-frequency modes are moderately damped. Regarding the primary mode (black line), for low wave--numbers, it follows the quadratic curve that governs the analytical solution. The curvature of this parabola depends linearly on the Peclet number ($Im(k)\sim -k^2/Pe$). As pointed out before, the amount of dissipation varies linearly with the inverse of the Peclet number, without affecting on the shape of the curves.

\subsubsection{Von Neumann analysis: effect of combined discretisation of viscous terms (viscosity $\mu$) and inviscid upwinding ($\lambda \text{diss}$)}\label{subsubsec:results_LESandUpwind}

\begin{figure}
	\centering
	\subfigure[Dispersion error.]{\includegraphics[scale=0.23]{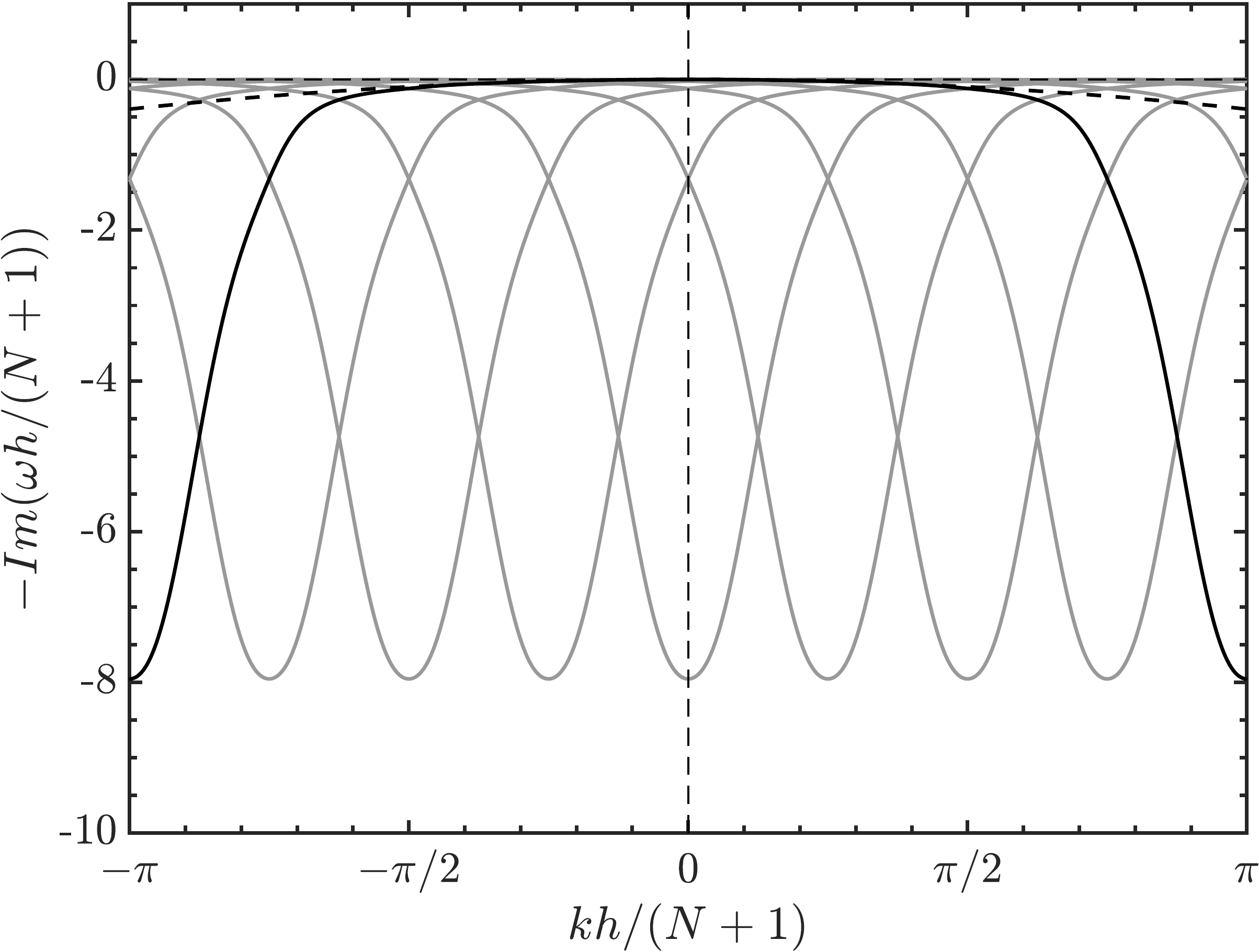}\label{fig:LambdaBR1_dispersion}}
	\subfigure[Dissipation error.]{\includegraphics[scale=0.23]{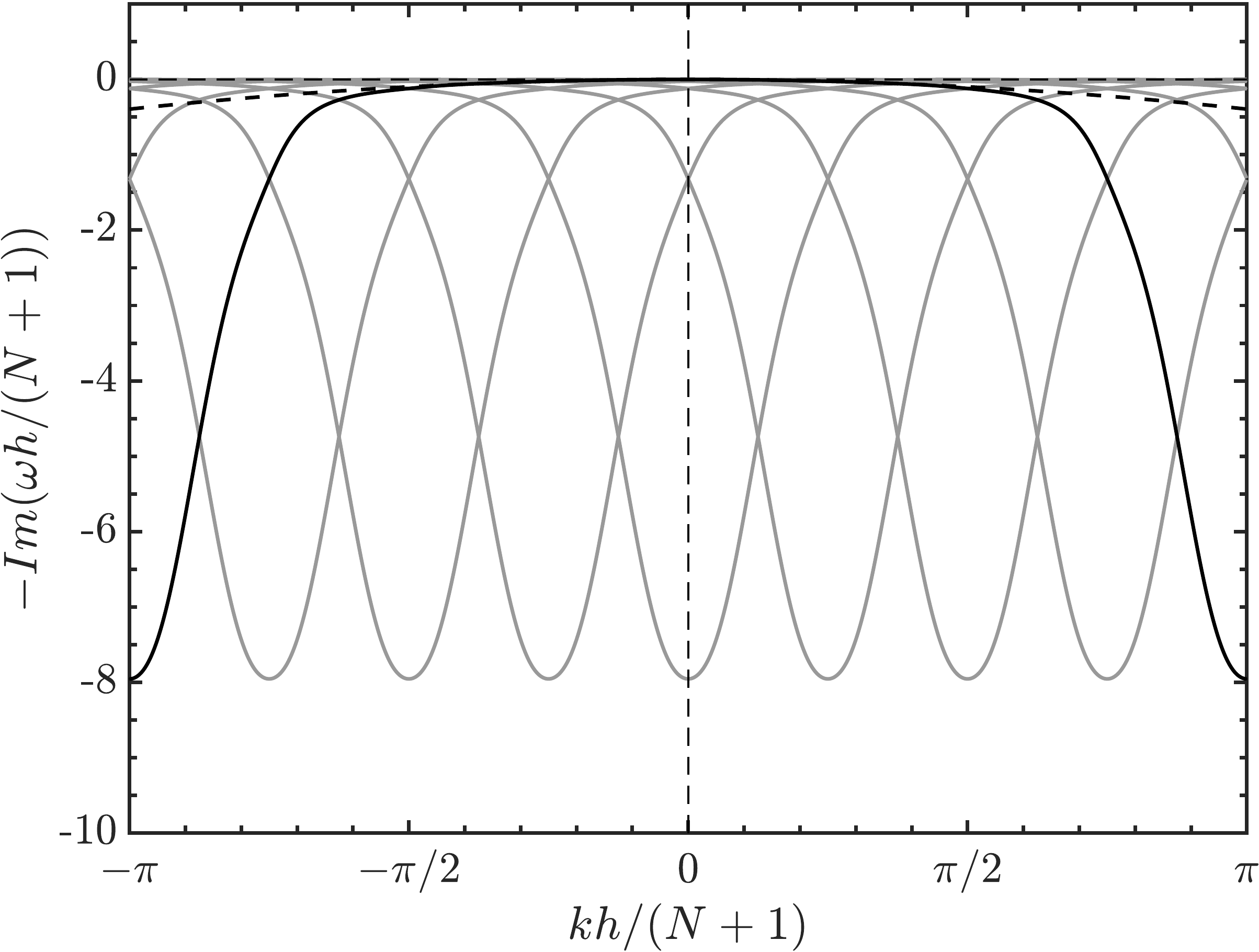}\label{fig:LambdaBR1_dissipation}}
	\caption{Dispersion and dissipation errors obtained with artificial viscosity (BR1 scheme) and upwind fluxes ($\lambda=1$). In this case, the Peclet number is $\text{Pe}=100$.}	
	\label{fig:LambdaBR1_dispdiss}	
\end{figure}

\begin{figure}
	\centering
	\subfigure[Full view.]{\includegraphics[scale=0.23]{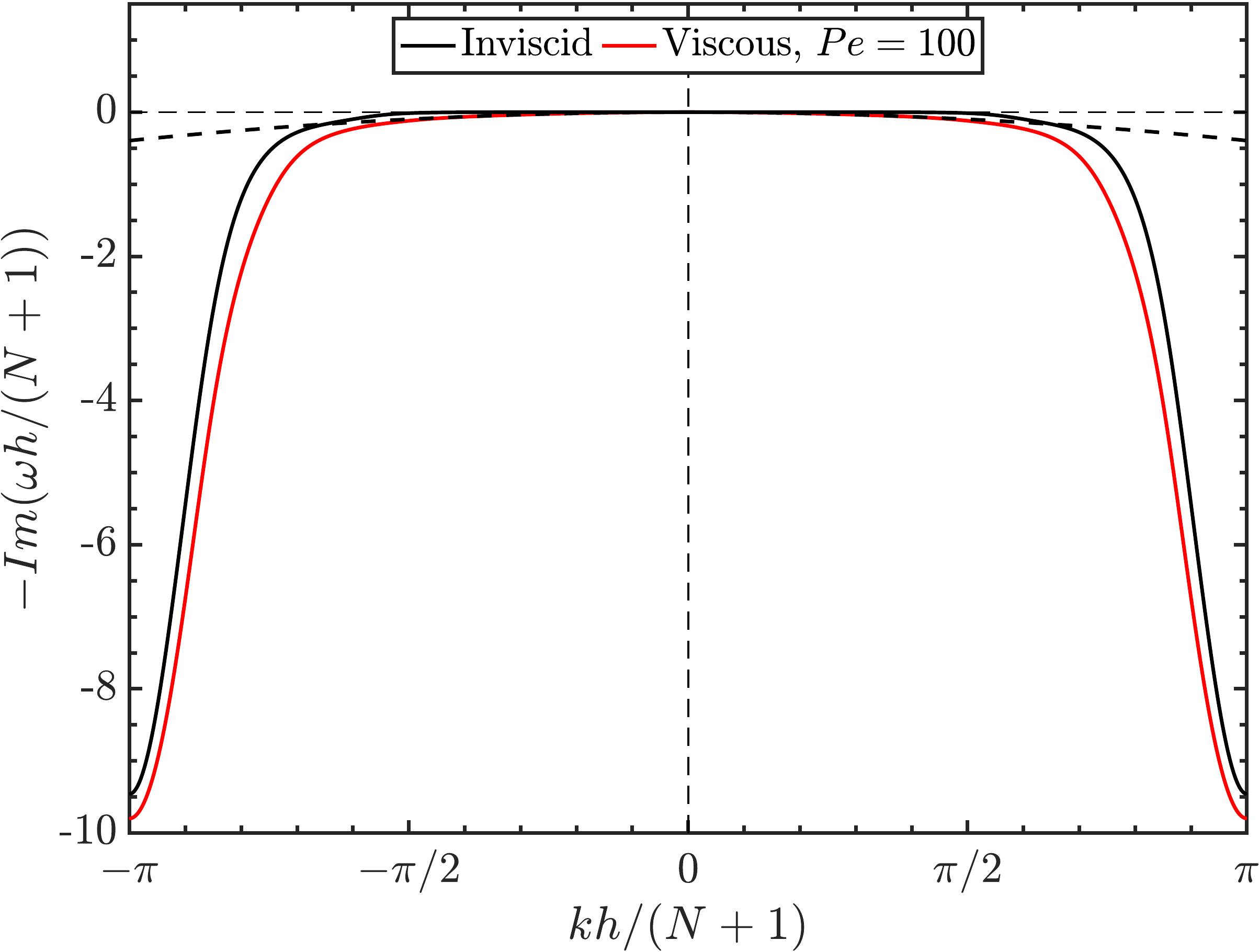}\label{fig:LambdaBR1_comparison}}
	\subfigure[Detailed view.]{\includegraphics[scale=0.23]{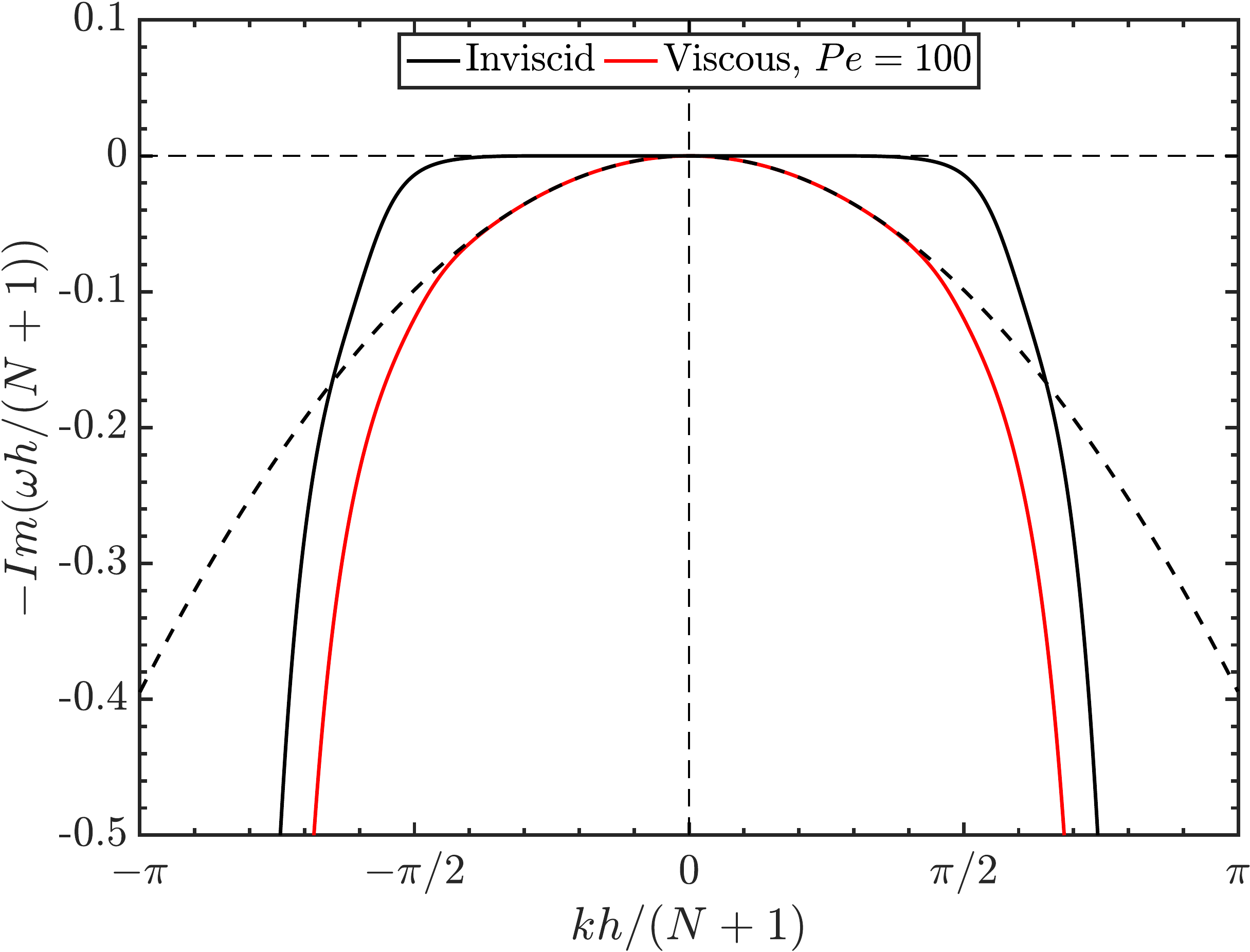}\label{fig:LambdaBR1_detail}}	
	\caption{Comparison of the inviscid and viscous primary mode when $\lambda=1$ and $Pe=100$. The dashed line shows the theoretical dissipation of second order derivatives, $k^2/Pe$.}
	\label{fig:LambdaBR1_comp}
\end{figure}

We now undertake an assessment on the effect in the dispersion and dissipation errors of combining interface stabilisation $\lambda \text{diss}(u,\phi)$ and in the advection--diffusion equation (controlled with the Peclet number, $Pe$). This also can be understood as increasing the penalty parameter, $\sigma$, in the Interior Penalty (IP) formulation \cite{2017:Ferrer}. We consider upwind ($\lambda=1$) for the inviscid flux discretisation, whilst maintaining the BR1 scheme for the artificial viscosity, now with Peclet number $Pe=100$ in an attempt to model flows at higher Reynolds numbers. Von Neumann results are shown in Figure \ref{fig:LambdaBR1_dispdiss}. Dispersion errors are depicted in Figure \ref{fig:LambdaBR1_dispersion}, where we notice that, as in the inviscid scheme, all secondary modes are replications of the primary mode. The dissipation is represented in Figure \ref{fig:LambdaBR1_dissipation}, where the dashed line represents the analytical solution \eqref{eq:AnalyticalSolutionWithPeclet} dissipation (for $Pe=100$). The comparison between these results and the inviscid solution primary mode dissipation error (in Figure \ref{fig:fluxes_1000_dissipation}) is detailed in Figure \ref{fig:LambdaBR1_comp}. In Figure \ref{fig:LambdaBR1_comparison} we find that at high wave--numbers, the impact of viscosity can be regarded as negligible compared to the amount of dissipation introduced by the Riemann solver. However, the detailed view, Figure \ref{fig:LambdaBR1_detail}, shows that the impact of the discretisation of the viscous terms is concentrated at small wave--numbers. In this region, the dissipation introduced by the upwind Riemann solver is negligible compared to that introduced by the discretisation of the viscous terms. Therefore, the overall dissipation in low wave--numbers is dominated by the theoretical bound ($k^2/Pe$) dictated by the viscous discretisation.

To conclude, the results discussed so far indicate that fluxes type stabilisation is mandatory to control dissipation at high wave--numbers, whilst the viscous terms govern dissipation at low and medium wave--numbers. This fluxes stabilisation can be introduced by both inviscid Riemann solvers $F^\star_e$ or other second order derivatives discretizations (e.g. IP).

\subsubsection{Navier--Stokes TGV problem}
\begin{figure}
\centering
\includegraphics[scale=0.3]{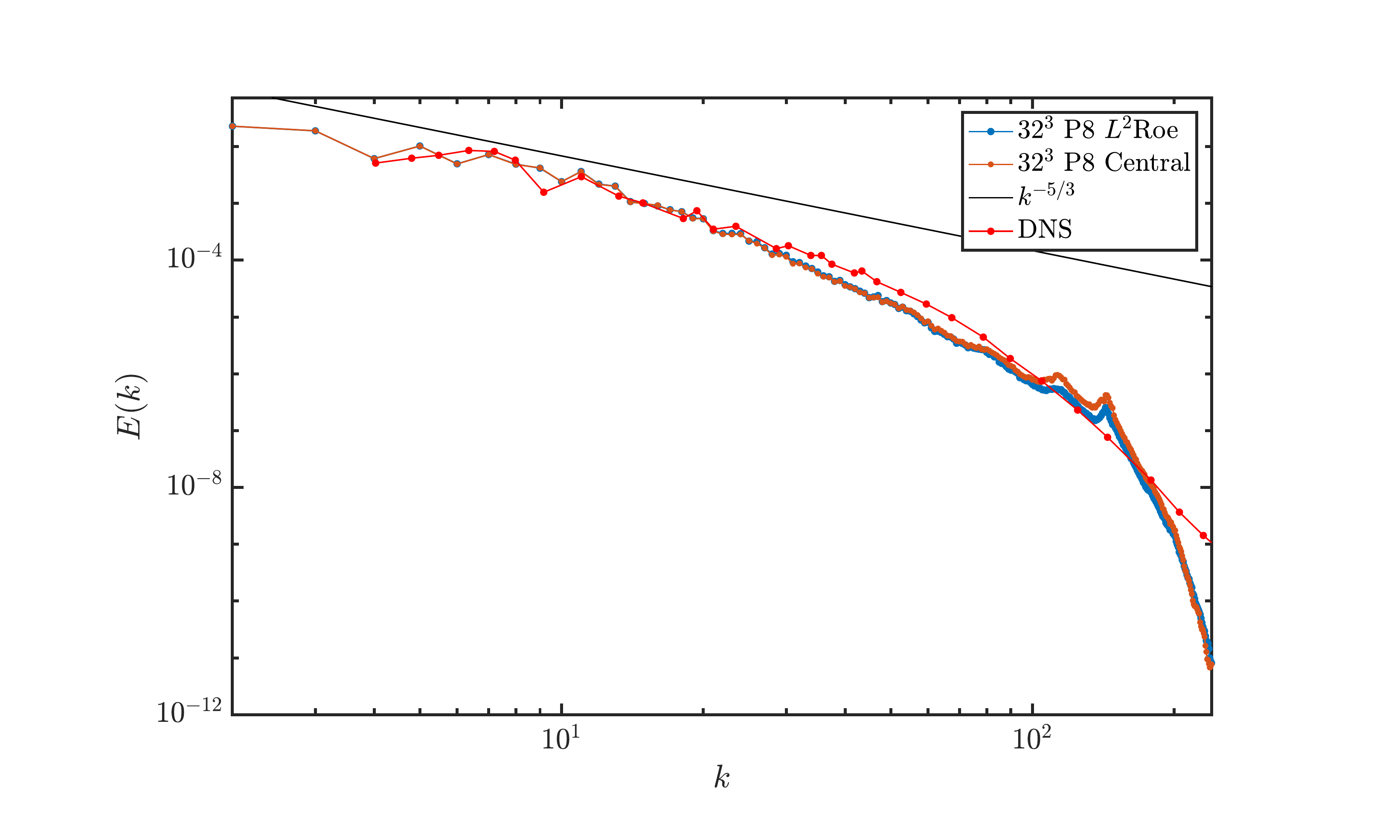}
\caption{Kinetic energy spectra obtained for the viscous TGV problem with $Re=1600$ and Mach number $M_0=0.1$.  The DG configuration consists in a linear mesh with $32^3$ elements, and polynomial order $N=8$. The figure shows the results obtained with central fluxes ($\lambda=0$), and with the low dissipation Roe Riemann solver presented in \cite{2016:Obwald}.}
\label{fig:BR1withRiemann}
\end{figure}

We consider the viscous Taylor--Green vortex problem with $Re=1600$ and Mach number $M_0=0.1$. The viscous version of the TGV is required to study viscous discretisations. We use $32^3$ elements with polynomial order $N=8$, which is a configuration close to the Direct Numerical Simulation (DNS). The kinetic energy spectra obtained in $t=8$ is represented in Figure \ref{fig:BR1withRiemann}, where we show the results obtained with central fluxes (pure viscous discretisation) and with the low dissipation Roe Riemann solver ($\lambda\sim0.1$, in equation \eqref{eq:RiemannStructureWithLambda}) introduced in \cite{2016:Obwald}. We find that without interface stabilisation, even with this detailed resolution level, there is a small accumulation of energy at high wave--numbers. This was previously noticed in \cite{2017:Flad}, where the authors explored the capabilities of different stabilisation techniques with an isotropic turbulent decay. It should be mentioned that the effect of this missing dissipation is more pronounced in the test case found in \cite{2017:Flad} that in the test case presented here. To obtain satisfactory results for coarser resolutions, numerical dissipation that affects low and medium wave--numbers will be required.

\subsection{Third dissipation technique: artificial viscosity (LES models and spectral vanishing viscosity)}\label{subsec:results_SVV}
Heretofore we have studied upwind Riemann solvers and the discretisation of viscous terms. The former has been found useful to add dissipation in high wave--numbers, whilst the latter adds dissipation in low and medium wave--numbers. In an attempt to obtain a more versatile scheme and a greater degree of control of the numerical dissipation, we study artificial viscosity (LES models and Spectral Vanishing Viscosity, SVV, method), applied to discontinuous Galerkin schemes.

\subsubsection{Von Neumann analysis}

This section is devoted to the assessment of artificial viscosity (LES models and spectral vanishing viscosity). We will perform von Neumann analysis to \eqref{eq:AdvEq} with the artificial viscous flux described in \eqref{eq:advectionSVV} using the power kernel written in \eqref{eq:SVVKernel}. It should be noticed that \eqref{eq:advectionLES} is recovered for $P_{SVV}=0$.
In Section \ref{subsubsec:results_LESandUpwind} it was highlighted the capability of the discretisation of the viscous terms to add dissipation in the low range wave--numbers (see Figure \ref{fig:LambdaBR1_detail} for details). This is mainly the effect of a simple LES model in the discretisation. However, we find that these methods tend to introduce non--vanishing dissipation even with smooth flows, as reported in \cite{2017:Fernandez}. The SVV is specially useful to overcome this drawback, since it filters--out low frequency modes and prevents dissipation of smooth solutions.

\begin{figure}
\centering
\subfigure[Full view]{\includegraphics[scale=0.25]{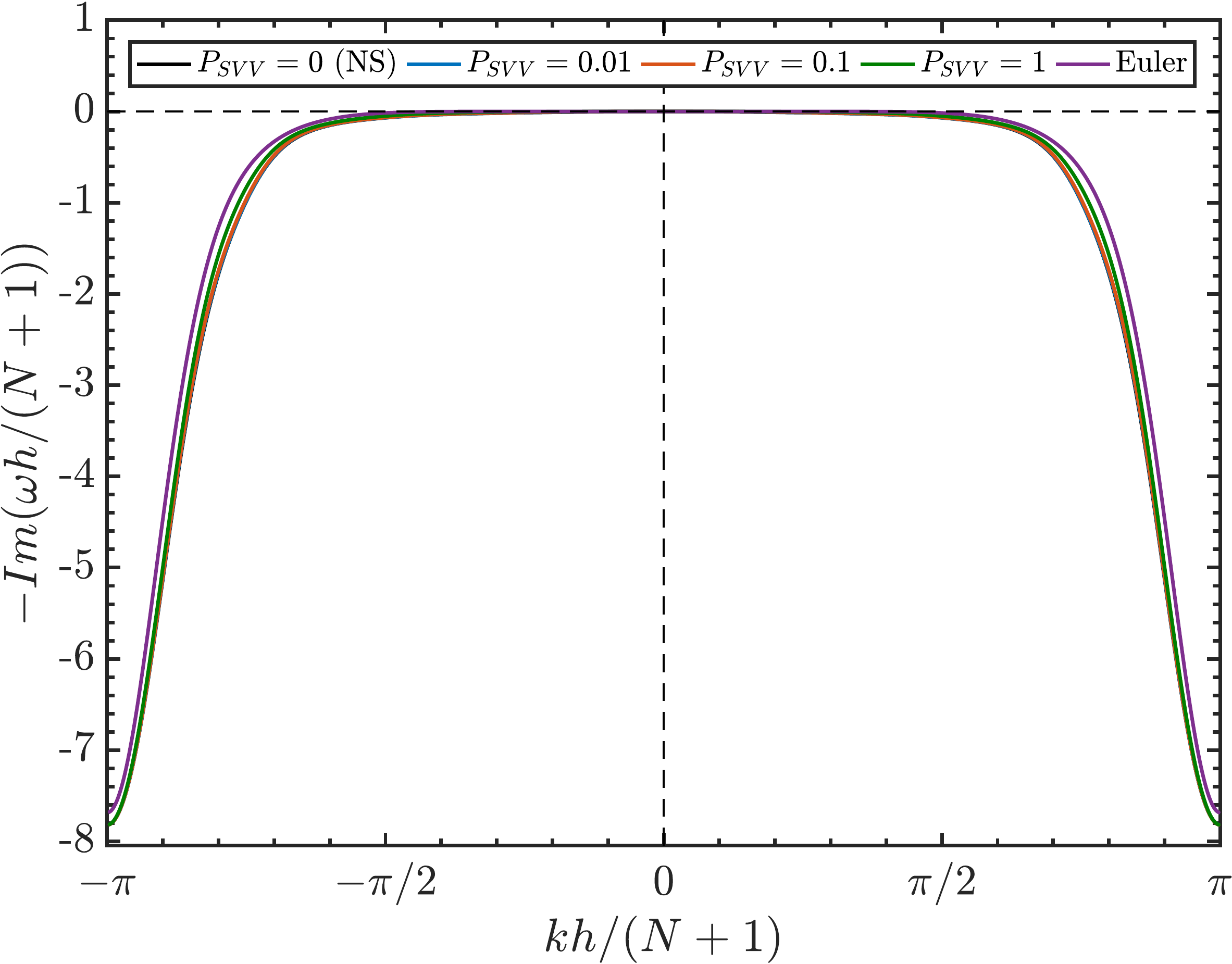}\label{fig:SVV}}
\subfigure[Detailed view]{\includegraphics[scale=0.25]{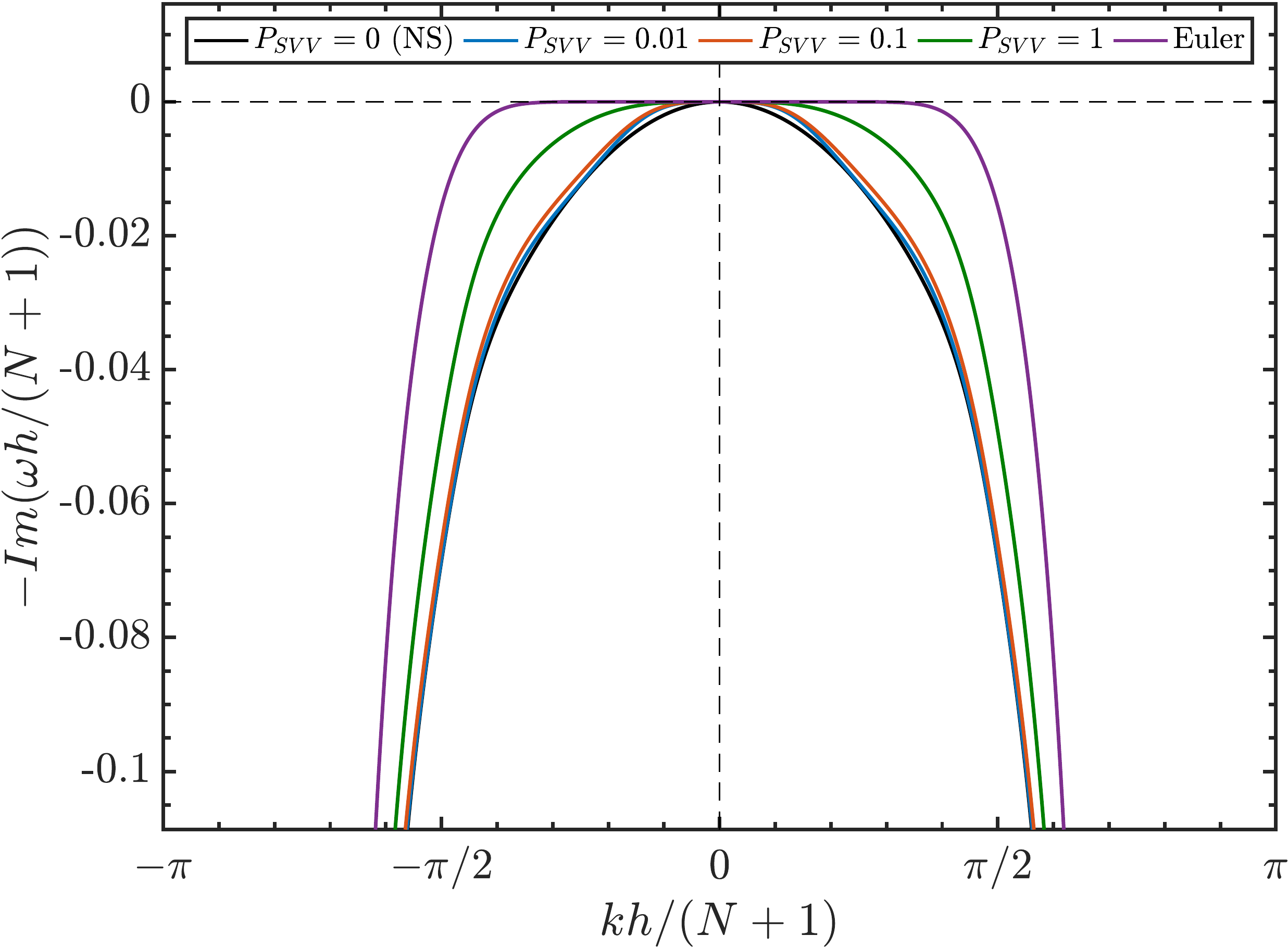}\label{fig:SVVdetail}}
\caption{Von Neumann dissipation curves using SVV with upwind ($\lambda=1$). Particularly, the effect of the SVV kernel power, $P_{SVV}$, is studied, where we have also included the particular cases with $P_{SVV}=0$ (standard second order derivative), and $P_{SVV}=\infty$ (inviscid, or Euler).}
\end{figure}

We depict von Neumann dissipation curves obtained with $\mu_{SVV}=0.005$ ($Pe_{SVV}=aL/mu_{SVV}=200$) and varying $P_{SVV}$ in Figure \ref{fig:SVV}, and its detailed view in Figure \ref{fig:SVVdetail}. In the default view  of Figure \ref{fig:SVV}, only a subtle difference is recognised since the effect of the SVV is concentrated in low and medium wave--numbers. Thus, to understand the effectiveness of the SVV, we check the dissipation curve concentrated on low and medium wave--numbers (Figure \ref{fig:SVVdetail}). We show that the SVV, with the power kernel \eqref{eq:SVVKernel}, is capable to control the shape of the dissipation curve between a pure viscous discretisation ($P_{SVV}=0$, labelled as NS) and the inviscid ($P_{SVV}>>1$, labelled as Euler), thus allowing a precise control of the numerical dissipation introduced. The impact of the SVV in the dissipation may seem negligible since the overall difference between all curves is small. Nevertheless, we will show by means of numerical experiments that these small differences are capable to effectively control the energy drain in under--resolved turbulent flow simulations.

\subsubsection{Navier--Stokes problem}

The linear von Neumann analysis shows that the SVV method allows to shape the dissipation curve between those obtained with upwind Riemann solvers and constant viscosities. Therefore, the SVV introduces dissipation mostly in low and medium wave--numbers. To test the SVV capabilities in Euler equations, we consider the inviscid Taylor--Green vortex problem described in \ref{subsec:theory_TGV}. In this study we use a coarser Cartesian $4^3$ mesh with $N=8$ to test the SVV in a severely under--resolved configuration. In this section, we have adopted the low dissipation Roe \cite{2016:Obwald} as Riemann solver (equivalent to $\lambda=0.1$ with $M_0=0.1$) as it was shown before (see Fig. \ref{fig:numExp:lambdaStudy:spectra}) that the standard Roe ($\lambda=1$) results in over-dissipated solutions for this problem. The SVV viscosity is set to $\mu_{SVV}=0.005$ (equivalent to a Reynolds number $Re_{SVV}=200$), and we study the effect of the kernel power coefficient by evaluating three values: $P_{SVV}=0.1, 1, $ and $10$. We have chosen this value for the SVV viscosity since in the NSE ($P_{SVV}=0$) it provides an excess of dissipation for the given mesh. Therefore, the dissipation excess is removed with the SVV filtering (controlled with $P_{SVV}$), instead of decreasing the viscosity $\mu_{SVV}$ to a more precise value.

\begin{figure}
\centering
\includegraphics[scale=0.23]{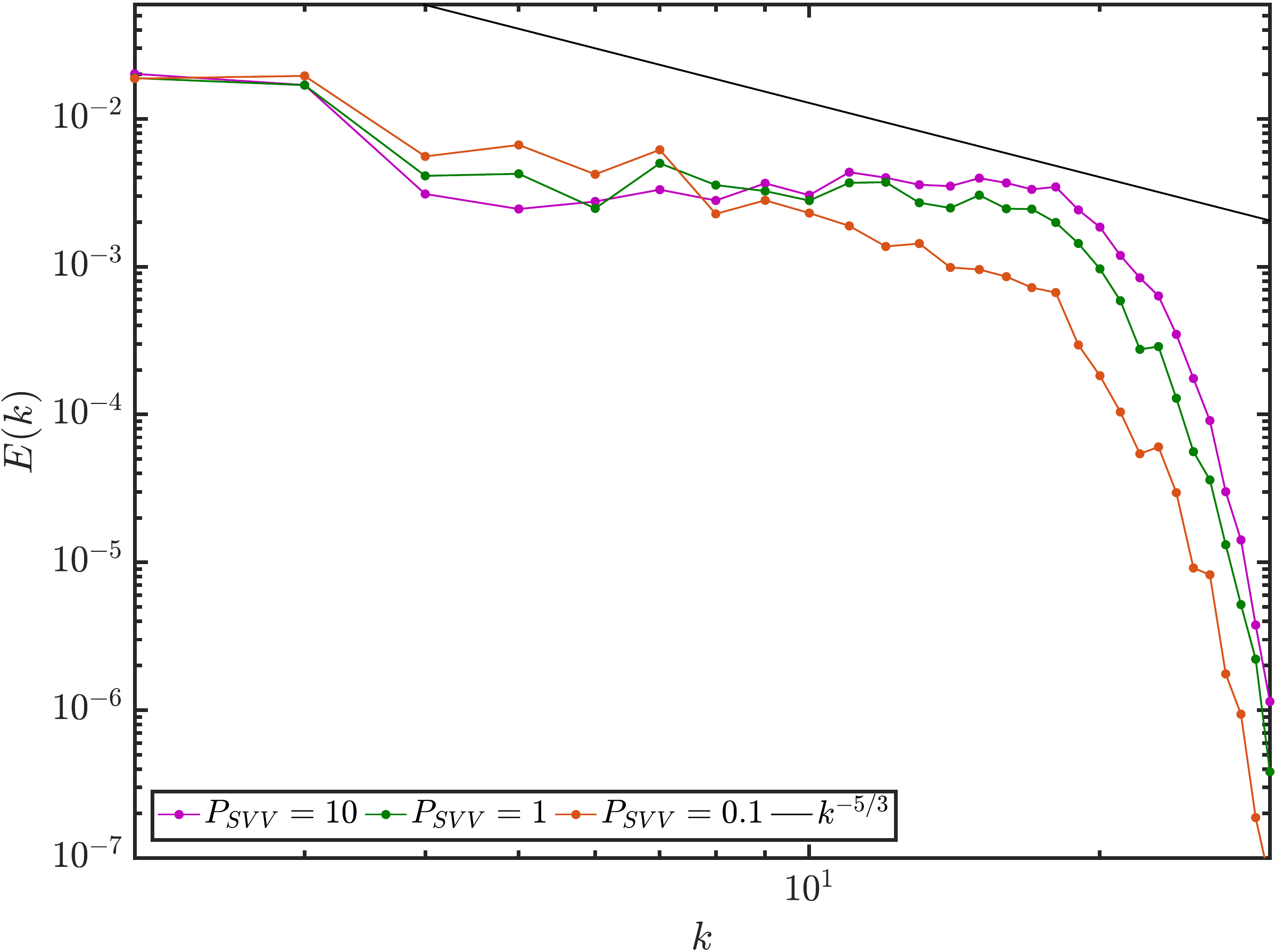}
\caption{Inviscid Taylor--Green vortex problem (with Mach number $M_0=0.1$) with $4^3$ elements and polynomial order $N=8$. We have used the SVV operator as defined in \eqref{eq:SVVNS}, where we set the SVV viscosity to $\mu_{SVV}=0.005$ (equivalent to a $Re_{SVV}=200$), and we vary the kernel power $P_{SVV}$.  Note that the major difference between the standard Navier--Stokes discretisation ($P_{SVV}=0$) and the $P_{SVV}=0.1$ scheme is dominant on low wave--numbers.}
\label{fig:numExps:TGV4P8SVV}
\end{figure}

The kinetic energy spectra in $t=8$ is depicted in Figure \ref{fig:numExps:TGV4P8SVV} for the three SVV kernel power $P_{SVV}$ values considered. We find that the energy decay in high wave--numbers is similar for all simulations, thus supporting von Neumann results in Figure \ref{fig:SVV}. In low and medium wave--numbers, we confirm the effectiveness of the SVV to shape the dissipation curve and adjust the energy spectra. We show that the value $P_{SVV}=0.1$ yields a satisfactory result according to Figure \ref{fig:numExps:TGV4P8SVV}, where the energy decays approximately following Kolmogorov's theoretical $k^{-5/3}$ rate \cite{pope2001turbulent}, and then is dissipated for high wave--numbers without producing energy accumulation. We find that both $P_{SVV}=1$ and $P_{SVV}=10$ present a lack of dissipation at medium wave--numbers, providing an implicit LES configuration with unsatisfactory results.

Hence, the resulting scheme combining upwind Riemann solvers and SVV is versatile, but requires an appropriate estimation of the SVV viscosity $\mu_{SVV}$ and the kernel power $P_{SVV}$. Besides, a proper LES model (variable value of the artificial viscosity) has not been tested yet. We address these issues in the following section.

\section{Design of a Smagorisnky-SVV scheme}\label{sec:LESSVV_hybrid}
In Section \ref{subsec:results_upwindRS} we have assessed the capabilities of upwind Riemann solvers to provide numerical dissipation in high wave--numbers, the capability of second order derivatives to introduce numerical dissipation at low and medium wave--numbers, and the potential of the SVV to shape the dissipation curve between pure elliptical discretisations ($P_{SVV}=0$) and the inviscid equation. The latter has been found to effectively introduce dissipation at medium wave--numbers and (if required) also at low wave--numbers. In this section, we combine the strategies presented in the previous section to construct a scheme capable to provide accurate solutions in turbulent under--resolved flows. Hence, we study the effectiveness of a scheme combining a low dissipation Riemann solver to damp high wave--numbers and a SVV to dissipate medium wave--numbers to solve the inviscid TGV problem.

One of the drawbacks of the SVV is its requirement to estimate two parameters with a remarkable impact on the final solution, (see Figure \ref{fig:numExps:TGV4P8SVV}). In an attempt to automatise the parameter selection, we set the SVV viscosity in \eqref{eq::SVVNS_tau} to that specified by a standard Smagorinsky LES turbulence model as suggested in \cite{2000:Karamanos} and implemented in \cite{2002:Kirby} in the context of continuous Galerkin methods:

\begin{equation}
\mu_{SVV} = C_S^2\Delta^2 |\boldsymbol{S}|,\label{eq:muSVVLES}
\end{equation}
being $S_{ij}=\frac{1}{2}(\partial_j v_{i} + \partial_i v_j)$ the strain tensor, and such that only the $P_{SVV}$ parameter remains free.

\begin{figure}
\centering
\subfigure[Energy spectra with polynomial order $N=4$]{\includegraphics[scale=0.23]{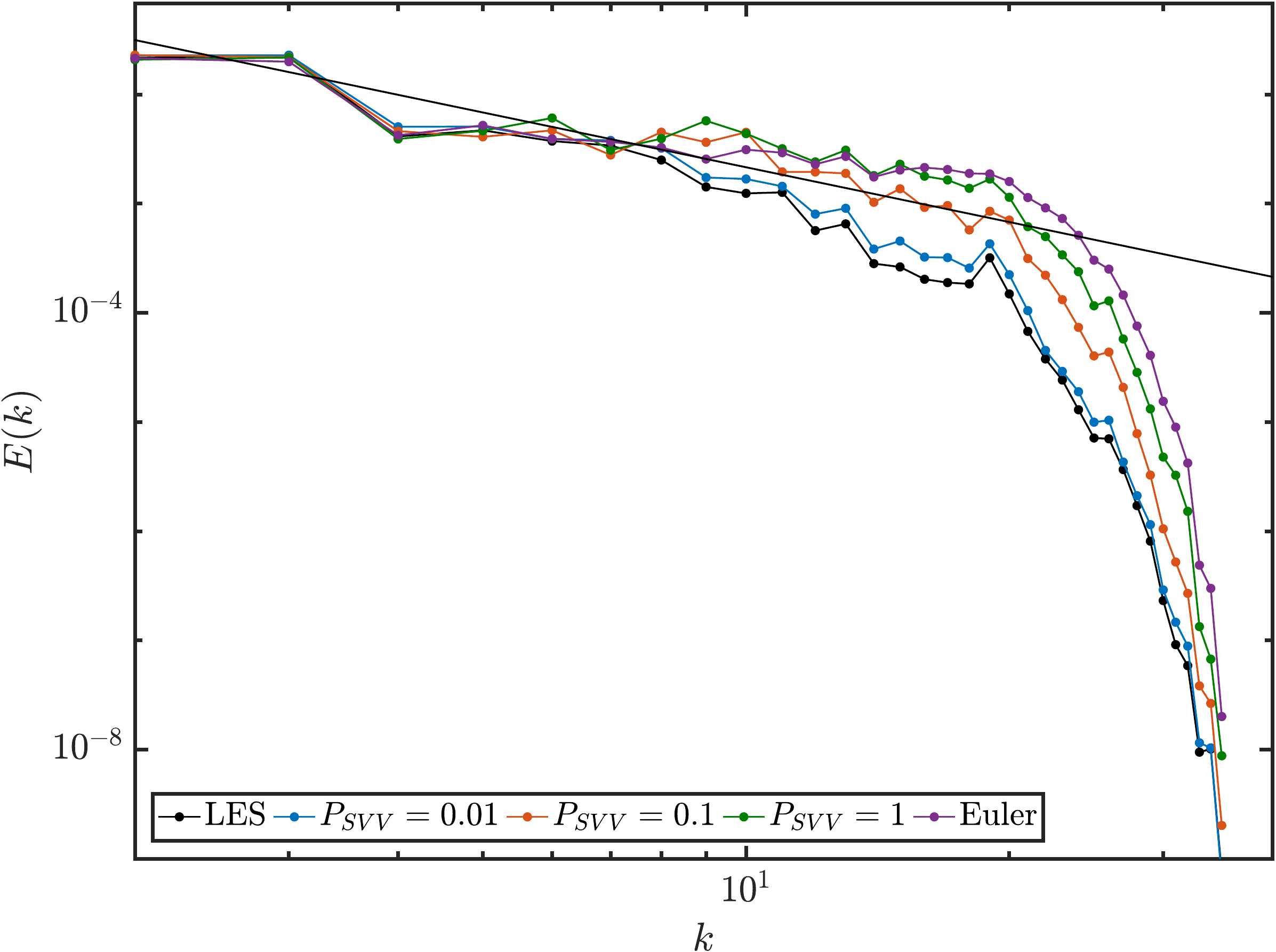}\label{fig:numExps:LES+SVV8P4spectra}}
\subfigure[Energy spectra with polynomial order $N=8$]{\includegraphics[scale=0.23]{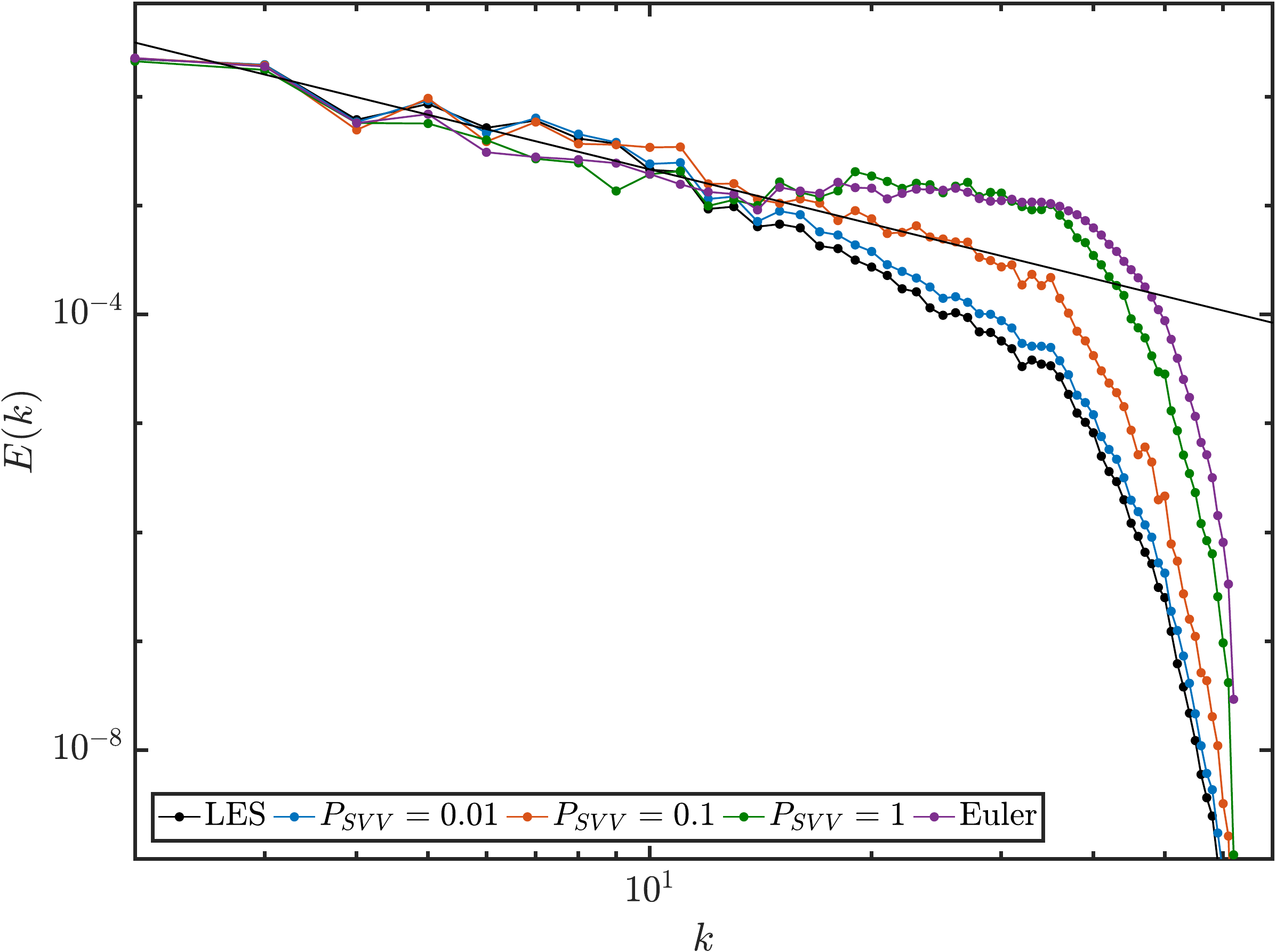}\label{fig:numExps:LES+SVV8P8spectra}}
\caption{Kinetic energy spectra in $t=8$ obtained with the proposed Smagorinsky--SVV strategy for two polynomial orders. Both cases were computed using a $8^3$ Cartesian mesh, solving the inviscid Taylor--Green vortex problem with Mach number $M_0=0.1$. Different values of the SVV kernel power $P_{SVV}$ were studied.}
\label{fig:numExps:LES+SVVspectra}
\end{figure}

\begin{figure}
\centering
\subfigure[Numerical viscosity with polynomial order $N=4$]{\includegraphics[scale=0.23]{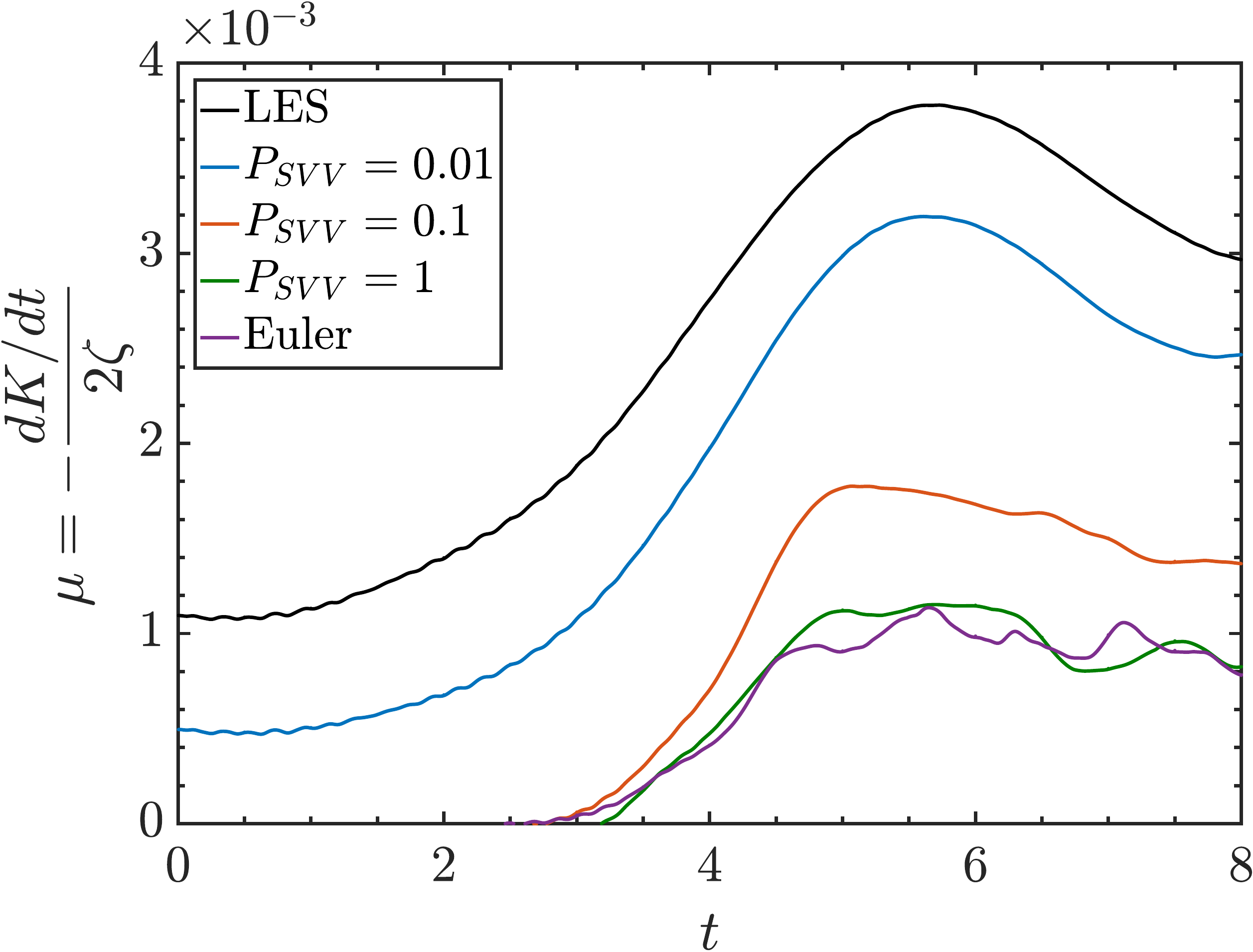}\label{fig:numExps:LES+SVV8P4visc}}
\subfigure[Numerical viscosity with polynomial order $N=8$]{\includegraphics[scale=0.23]{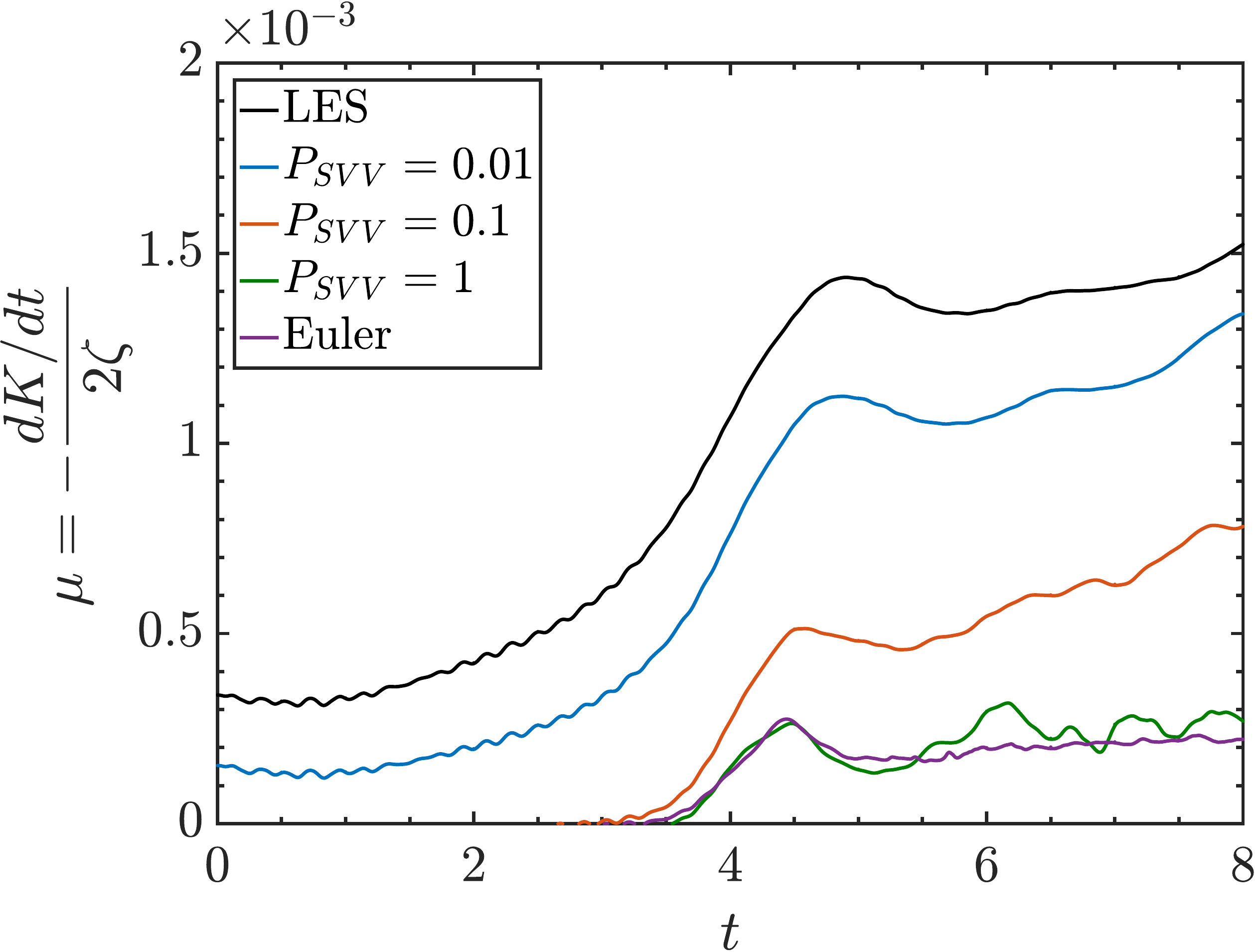}\label{fig:numExps:LES+SVV8P8visc}}
\caption{Numerical viscosity introduced by the Smagorinsky--SVV strategy. Both cases were computed using a $8^3$ Cartesian mesh. Different values of the SVV kernel power $P_{SVV}$ were studied. Two effects are regarded: the capability of the SVV to remove the dissipation in the laminar region of the TGV problem ($t<3$), and the lack of dissipation presented by the SVV with high kernel power coefficients.}
\label{fig:numExps:LES+SVVvisc}
\end{figure}

\begin{figure}
\centering
\includegraphics[scale=0.3]{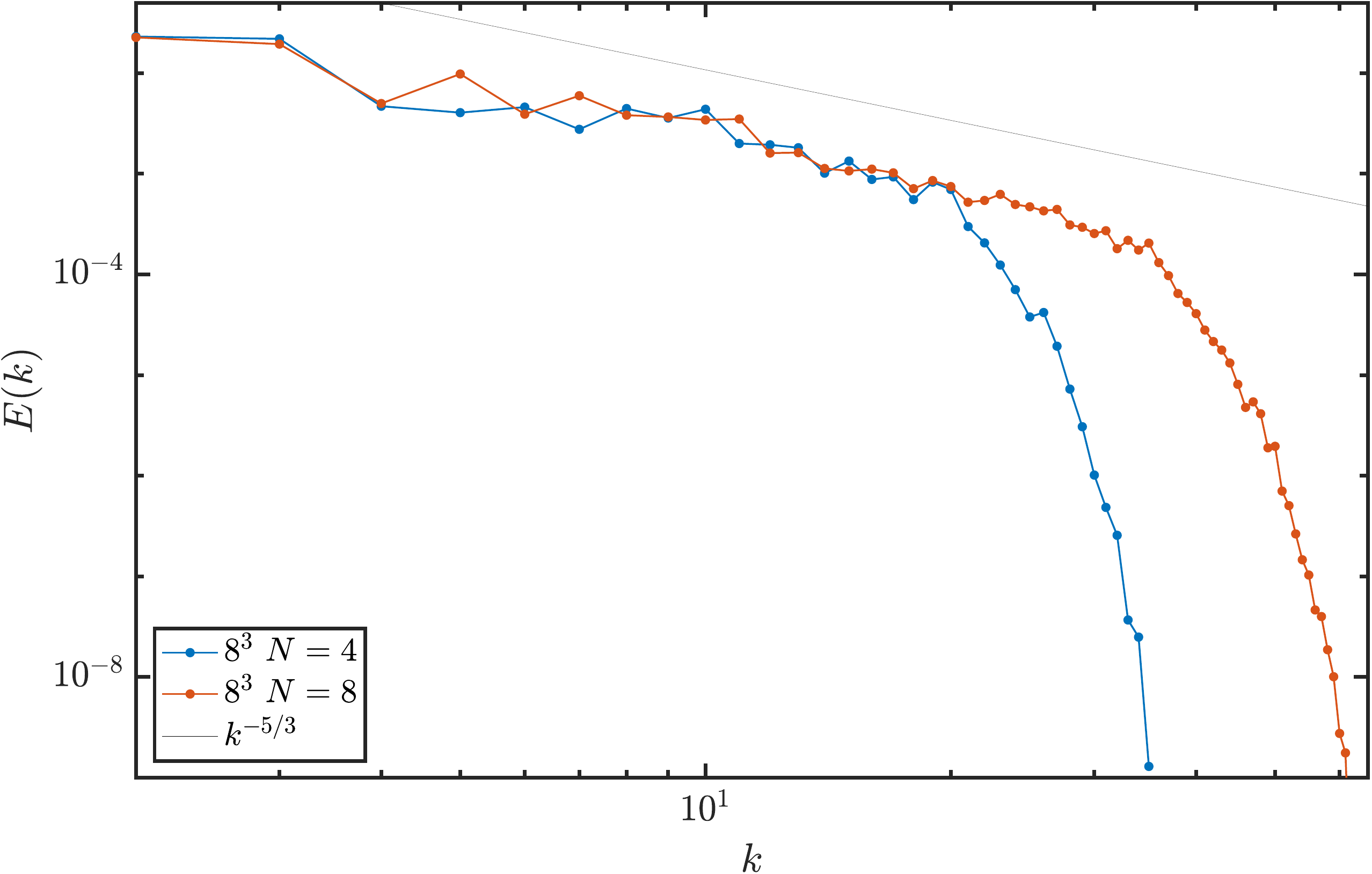}\label{fig:numExps:LES+SVVPsvv01spectra}
\caption{Inviscid Taylor--Green vortex problem ($M_0=0.1$) kinetic energy spectra in $t=8$ obtained with the Smagorinsky--SVV strategy using two different polynomial orders. Both cases were computed using a $8^3$ Cartesian mesh, and the SVV kernel power is $P_{SVV}=0.1$.}
\label{fig:numExps:LES+SVVPsvv01}
\end{figure}

We consider the inviscid Taylor--Green vortex problem in a Cartesian mesh with $8^3$ elements and two polynomial orders, $N=4$ and $N=8$, and we maintain the low dissipation Roe as Riemann solver. Both cases use the SVV method, with the viscosity computed as described in \eqref{eq:muSVVLES}. Figure \ref{fig:numExps:LES+SVVspectra} depicts the kinetic energy spectra varying the SVV kernel power $P_{SVV}$, and the standard Smagorinsky model (without SVV) for comparison. We show that the resulting spectra, as in Figure  \ref{fig:numExps:TGV4P8SVV} shows high dependence to the SVV power kernel, $P_{SVV}$. Precisely, for values $P_{SVV}>0.1$ the dissipation provided by the method is not enough to control the accumulation of energy at high wave--numbers. We show that the main difference with the standard Smagorinsky model is that the latter presents an excessive dissipation when the flow is laminar, hence, decreasing the overall energy at further times. This was also observed by \cite{2017:Fernandez}. This is naturally avoided with the SVV technique, as it filters--out the laminar (smooth) energy components of the dissipation, introducing only the high frequency fluctuations. For completeness, we represent the numerical viscosity introduced by both Smagorinsky (LES) and Smagorinsky--SVV in Figure \ref{fig:numExps:LES+SVVvisc}. First, the laminar region (the region without subgrid--scales in the flow, $t<3$) presents non--negligible numerical dissipation when using the standard Smagorinsky model, whilst this dissipation vanishes when considering the Smagorinsky--SVV approach. Second, we confirm the lack of dissipation presented by the SVV method with excessive kernel power $P_{SVV}$. As a conclusion, we find the value $P_{SVV}=0.1$ appropriate for this test case. Lastly, we show the energy spectra for this configuration in Figure \ref{fig:numExps:LES+SVVPsvv01}, showing that the turbulence model is not altered by the polynomial order. In summary, the low dissipation Roe Riemann solver, the Smagorinsky model (which has been recently adopted in \cite{2017:Flad} as LES model in the context of DG solvers) as an estimation of the $\mu_{SVV}$ viscosity, and $P_{SVV}=0.1$ provides a robust method to simulate isotropic turbulence.

\section{Conclusions}\label{sec:conclusions}

In this paper we have studied the dissipation introduced by three strategies: upwind Riemann solvers, discretisation of viscous terms and artificial viscosity (LES models and spectral vanishing viscosity). We have performed a linear analysis using von Neumann method to discover that the dissipation introduced by upwind Riemann solvers  is not linear with parameter $\lambda$, which penalises the inter--element solution jumps. Instead, the dissipation introduced increases until a critical value is reached, to then decrease as the discretisation tends to that of a conforming (i.e. continuous Galerkin) method. We found that an upwind Riemann solver is required to maintain the dissipation introduced at high wave--numbers.

Next, we studied the dissipation introduced by second order derivatives approximated with the BR1 method, to conclude that the dissipation at low and medium wave--numbers can be controlled with viscous discretisations, whilst Riemann solvers are entirely responsible of introducing dissipation in high wave--numbers. Lastly, we study the dissipation introduced by the spectral vanishing viscosity. We found that, with an appropriate filtering kernel, it is possible to achieve high control on the dissipation introduced in the low and medium wave--number region, maintaining low dissipation when the flow is laminar.

We compared the results obtained from von Neumann analyses of a linear equation, with simulations of the non--linear Euler/Navier--Stokes equations. To do so, we used the Taylor--Green Vortex (TGV) problem to test the behaviour of the different configurations in transitional/turbulent under--resolved flows. First, we found that the features of the different operators predicted by von Neumann analysis were also present in the Navier--Stokes equations. Second, we show that introducing the SVV method combined with a Smagorinsky model, and a low dissipation Riemann solver appears as the most robust and accurate technique to simulate isotropic under--resolved turbulence. This novel DG--SVV method is capable of maintaining low dissipation levels in laminar flows, whilst providing enough dissipation for increased Reynolds numbers.

\section*{Acknowledgments}
The authors acknowledge the computer resources and technical assistance provided by the Centro de Supercomputaci\'on y Visualizaci\'on de Madrid (CeSViMa).

\appendix
\section{Compressible Navier--Stokes formulation}\label{appendix::A}

The 3D Navier--Stokes equations can be compactly written as:
%
\begin{equation}
\vec{u}_t  + \nabla\cdot\vec{\boldsymbol{F}}_e = \nabla\cdot\vec{\boldsymbol{F}}_v,
\end{equation}
where $\vec{u}$ is the vector of conservative variables $\vec{u} = [ \rho , \rho v_1 , \rho v_2 , \rho v_3 , \rho e]^T$, $\vec{\boldsymbol{F}}_e$ are the inviscid, or Euler equations fluxes:

\begin{equation}
\vec{\boldsymbol{F}_e} = \left[\begin{array}{ccc} \rho v_1 & \rho v_2 & \rho u_3 \\
                                                                                \rho v_1^2 + p & \rho v_1v_2 & \rho v_1v_3 \\
                                                                                	\rho v_1v_2 & \rho v_2^2 + p & \rho v_2v_3 \\
                                                                                	\rho v_1v_3 & \rho v_2v_3 & \rho v_3^2 + p \\
                                                                                	\rho v_1 H & \rho v_2 H & \rho v_3 H
\end{array}\right],
\end{equation}
where $\rho$, $e$, $H$ and $p$ are the density, total energy, total enthalpy, and pressure respectively, and  $\vec{v}=[v_1,v_2,v_3]^T$ is the velocity. Additionally, ${\boldsymbol{F}}_v$ defines the viscous fluxes:

\begin{equation}
\vec{\boldsymbol{F}_v}(\mu,\vec{v},\nabla\vec{v}) = \left[\begin{array}{ccc}0 & 0 & 0\\
 \tau_{xx} & \tau_{xy} & \tau_{xz} \\
 \tau_{yx} & \tau_{yy} & \tau_{yz} \\
 \tau_{zx} & \tau_{zy} & \tau_{zz} \\
 \sum_{j=1}^3 v_j\tau_{1j} + \kappa T_x& \sum_{j=1}^3 v_j\tau_{2j} + \kappa T_y& \sum_{j=1}^3 v_j\tau_{3j} + \kappa T_z
\end{array}\right],
\label{eq:viscousfluxes}
\end{equation}
where $\kappa$ is the thermal conductivity, $T_x, T_y$ and $T_z$ denote the gradients of temperature and the stress tensor $\boldsymbol{\tau}$ is defined as $\boldsymbol{\tau} = \mu(\nabla \vec{v} + (\nabla \vec{v})^T) - 2/3\mu \boldsymbol{I}\nabla\cdot\vec{v}$, with $\mu$ the dynamic viscosity, and $\boldsymbol{I}$ is the three-dimensional identity matrix.

  We discretise Euler equations using the novel nodal DG split--formulation derived in \cite{2016:gassner} (precisely, we use Pirozzoli split formulation), whilst for viscous fluxes we use the Bassi-Rebay 1 (BR1) scheme \cite{1997:Bassy}. In all cases, the Mach number is kept to 0.1 such that compressible effects are negligible.

\section{Taylor--Green vortex problem}\label{subsec:theory_TGV}

Numerical experiments are performed to evaluate the validity of von Neumann analysis assessments in the more general case of the NSE. To test the under--resolved capabilities of the strategies studied with von Neumann analysis, we will solve the Taylor--Green Vortex (TGV) problem \cite{1937:Taylor}. The TGV problem has been widely used to report the subgrid--scale modelling capabilities of iLES approaches and discretisations \cite{2016:gassner,2017:Moura}. In this paper, we assess von Neumann analysis truthfulness to estimate the dissipation introduced by the operators introduced in Sections \ref{subsubsec:upwindRiemannTheory}, \ref{subsubsec:theory_LES} and \ref{subsubsec:theory_SVV}.
The configuration of the TGV problem is a three dimensional periodic box $[-\pi,\pi]^3$ with the initial condition:

\begin{equation}
\begin{split}
\rho &= \rho_0,\\
v_1 &= V_0 \sin x\cos y \cos z,\\
v_2 &= -V_0\cos x \sin y \cos z,\\
v_3 &= 0,\\
p &= \frac{\rho_0 V_0^2}{\gamma M_0^2} + \frac{\rho_0 V_0^2}{16}(\cos 2x + \cos 2y)(\cos 2z + 2).\\
\end{split}
\label{eq:TGV}
\end{equation}
The Mach number is $M_0=0.1$ in all the simulations performed herein. The reported quantities to measure the simulations accuracy are the kinetic energy rate:

\begin{equation}
\epsilon = -\frac{dK}{dt} = -\frac{1}{|\Omega|}\frac{d}{dt}\int_{\Omega}\frac{1}{2}\rho V^2 d\boldsymbol{x},
\end{equation}
the enstrophy:

\begin{equation}
\zeta = \frac{1}{2|\Omega|}\int_{\Omega}(\nabla\times \boldsymbol{v})^2 d\boldsymbol{x},
\end{equation}
the numerically introduced dissipation estimated with both $\epsilon$ and $\zeta$ \cite{2013:GassnerUnder}:

\begin{equation}
\mu \simeq  \frac{\epsilon}{2\zeta} = -\frac{dK/dt}{2\zeta},
\end{equation}
and the kinetic energy spectra, measured at a fixed time snapshot ($t=8$ to observe transitional flow, and $t=14$ to show the isotropic decay). In this paper we consider both the inviscid version of the Taylor--Green vortex problem, and the viscous Taylor--Green vortex problem with Reynolds number $Re=1600$.
Finally, all Navier--Stokes simulations are time-marched using a three stage Runge--Kutta scheme with a Courant--Friedrichs--Lewy (CFL) number of 0.4.

\bibliography{mybibfile}

\end{document}